\documentclass[11pt, a4paper]{article}

\usepackage{array}
\usepackage{multirow}

\makeatletter
\newcommand{\thickhline}{%
    \noalign {\ifnum 0=`}\fi \hrule height 1pt
    \futurelet \reserved@a \@xhline
}
\newcolumntype{"}{@{\hskip\tabcolsep\vrule width 1pt\hskip\tabcolsep}}
\makeatother

\usepackage{graphics}
\usepackage{epsfig}
\usepackage{subfigure}
\usepackage{amssymb}
\usepackage{amsthm}
\usepackage{mathrsfs}
\usepackage{hhline}
\usepackage{longtable}
\usepackage{amsmath}
\usepackage{array}
\usepackage{float}
\usepackage{picinpar}
\usepackage{multirow}
\usepackage{cite}
\usepackage{enumerate}
\usepackage{xcolor}
\usepackage[mathlines]{lineno}
\modulolinenumbers[2]
\newcommand*\patchAmsMathEnvironmentForLineno[1]{%
\expandafter\let\csname old#1\expandafter\endcsname\csname #1\endcsname  \expandafter\let\csname oldend#1\expandafter\endcsname\csname end#1\endcsname  \renewenvironment{#1}%
{\linenomath\csname old#1\endcsname}%
{\csname oldend#1\endcsname\endlinenomath}}%
\newcommand*\patchBothAmsMathEnvironmentsForLineno[1]{%
\patchAmsMathEnvironmentForLineno{#1}%
\patchAmsMathEnvironmentForLineno{#1*}}%
\AtBeginDocument{%
\patchBothAmsMathEnvironmentsForLineno{equation}%
\patchBothAmsMathEnvironmentsForLineno{align}%
\patchBothAmsMathEnvironmentsForLineno{flalign}%
\patchBothAmsMathEnvironmentsForLineno{alignat}%
\patchBothAmsMathEnvironmentsForLineno{gather}%
\patchBothAmsMathEnvironmentsForLineno{multline}%
}

\newtheoremstyle{rem}
  {10pt}          
  {10pt}  
  {\rm}  
  {}
  {\bf}  
  {: }    
  { }    
  {}     
\theoremstyle{rem}
\newtheorem{rem}{Remark}[section]
\newtheorem{example}{Example}[section]
\newtheorem{problem}{Problem}[section]

\newtheoremstyle{theorem}
  {10pt}          
  {10pt}  
  {\it}  
  {}
  {\bf}  
  {: }    
  { }    
  {}     
\theoremstyle{theorem}

\newtheorem{theorem}{Theorem}[section]
\newtheorem{lemma}[theorem]{Lemma}

\newtheorem{corollary}[theorem]{Corollary}
\newtheorem{conjecture}{Conjecture}[section]

\numberwithin{equation}{section}

\def\N{\mathbb{N}}
\newtheorem{defi}{Definition}[section]

\usepackage{graphicx}
\graphicspath{%
    {converted_graphics/}
    {/}
}

\def\ms{\medskip}
\def\nt{\noindent}
\begin{document}
\UseRawInputEncoding
\baselineskip18truept
\normalsize
\begin{center}
{\mathversion{bold}\Large \bf Every graph is local antimagic total and its application to local antimagic (total) chromatic numbers }

\bigskip
{\large  Gee-Choon Lau}\\

\medskip

\emph{Faculty of Computer \& Mathematical Sciences,}\\
\emph{Universiti Teknologi MARA (Segamat Campus),}\\
\emph{85000, Johor, Malaysia.}\\
\emph{geeclau@yahoo.com}\\

\end{center}

\begin{abstract}
A graph $G = (V, E)$ of order $p$ and size $q$ is said to be local antimagic if there exists a bijection $g:E(G) \to \{1,2,\ldots,q\}$ such that for any pair of adjacent vertices $u$ and $v$, $g^+(u)\ne g^+(v)$, where $g^+(u)=\sum_{uv\in E(G)} g(uv)$ is the induced vertex color of $u$ under $g$. We also say $G$ is local antimagic total if there exists a bijection $f: V\cup E \to\{1,2,\ldots ,p+q\}$ such that for any pair of adjacent vertices $u$ and $v$, $w(u)\not= w(v)$, where $w(u)= f(u) +\sum_{uv\in E(G)} f(uv)$ is the induced vertex weight of $u$ under $f$.  The local antimagic (and local antimagic total) chromatic number of $G$, denoted $\chi_{la}(G)$ (and $\chi_{lat}(G)$), is the minimum number of distinct induced vertex colors (and weights) over all local antimagic (and local antimagic total) labelings of $G$. We also say a local antimagic total labeling is local super antimagic total if $f(v)\in\{1,2,\ldots,p\}$ for each $v\in V(G)$. In [Proof of a local antimagic conjecture, {\it Discrete Math. Theor. Comp. Sc.}, {\bf 20(1)} (2018), \#18], the author proved that every connected graph of order at least 3 is local antimagic. Using this result, we provide a very short proof that every graph is local antimagic total. As an application, we showed that there exists close relationship between $\chi_{la}(G \vee K_1)$ and $\chi_{lat}(G)$. A sufficient condition is also given for the corresponding local super antimagic total labeling. Sharp bounds of $\chi_{lat}(G)$ and close relationships between $\chi_{lat}(G)$ and $\chi_{la}(G \vee K_1)$ are found. Bounds of $\chi_{lat}(G-e)$ in terms of $\chi_{lat}(G)$ for a graph $G$ with an edge $e$ deleted are also obtained. These relationships are used to determine the exact values of $\chi_{lat}(G)$, $\chi_{lat}(G-e)$, $\chi_{la}(G\vee K_1)$ and $\chi_{la}((G-e)\vee K_1)$ for many dense and sparse graphs $G$. The work of this paper also provides many open problems on $\chi_{lat}(G)$. We also conjecture that each graph $G$ of order at least 3 has $\chi_{lat}(G)\le \chi_{la}(G)$.\\


\end{abstract}

\section{Introduction}
Consider a $(p,q)$-graph $G(V,E)$ of order $p$ and size $q$. If $e$ is an edge of $G$, denote by $G-e$ the graph $G$ with the edge $e$ deleted. Let $g:E(G)\to [1,q]$ be a bijective edge labeling that induces a vertex labeling $g^+: V(G) \to \N$ such that $g^+(v) = \sum_{uv\in E(G)} g(uv)$. We say $g$ is a local antimagic labeling of $G$ if $g^+(u)\neq g^+(v)$ for each $uv\in E(G)$~\cite{Arumugam,Bensmail}. The minimum number of distinct vertex labels induces by $g$ is called the local antimagic chromatic number of $G$, denoted $\chi_{la}(G)$ \cite{Arumugam}. Clearly, $\chi_{la}(G)\ge \chi(G)$. Let $G \vee H$ be the join of $G$ and $H$ with vertex set $V(G)\cup V(H)$ and edge set $E(G)\cup E(H) \cup \{uv : u\in V(G), v\in V(H)\}$. Let $G+H$ be the disjoint union of $G$ and $H$ with vertex set $V(G)\cup V(H)$ and edge set $E(G) \cup E(H)$. For convenience, $nG$ is the disjoint union of $n\ge 1$ copies of $G$, and $nK_1 = O_n$. For positive integers $a$ and $b$, assume $[a,b]=\{a,a+1,\ldots,b\}$ and ${a\choose b}=0$ if $a < b$. 

\ms\nt  Let $f: V(G)\cup E(G) \to [1,p+q]$ be a bijective total labeling that induces a vertex labeling $w : V(G) \to \N$ such that $$w_f(v)=f(v) + \sum_{uv\in E(G)} f(uv)$$ and is called the weight of $v$ for each vertex $v \in V(G)$. We say $f$ is a local antimagic total labeling of $G$ (and $G$ is local antimagic total) if $w(u)\neq w(v)$ for each $uv\in E(G)$. Clearly, $w$ corresponds to a proper vertex coloring of $G$ if each vertex $v$ is assigned the color $w(v)$. Let $w(f)$ be the number of distinct vertex weights induces by $f$. The number $min\{w(f)\,|\,f \mbox{ is a local antimagic total labeling of } G\}$ is called the local antimagic total chromatic number of $G$, denoted $\chi_{lat}(G)$. Clearly, $\chi_{lat}(G)\ge \chi(G)$. It is well known that determining the chromatic number of a graph $G$ is NP-hard~\cite{Zukerman}. Thus, in general, it is very difficult to determine $\chi_{la}(G)$ and $\chi_{lat}(G)$. 

\ms\nt In~\cite{Haslegrave}, the author proved that every connected graph of order at least 3 is local antimagic. Using this result, we provide in Section~\ref{sec-pre} a very short proof that every graph is local antimagic total. As an application, we showed that there exists close relationship between $\chi_{la}(G \vee K_1)$ and $\chi_{lat}(G)$. A sufficient condition is also given for the corresponding local super antimagic total labeling.  Sharp bounds of $\chi_{lat}(G)$ and close relationships between $\chi_{lat}(G)$ and $\chi_{la}(G \vee K_1)$ are found. Bounds of $\chi_{lat}(G-e)$ in terms of $\chi_{lat}(G)$ are also obtained. These relationships are then used to determine the exact value of $\chi_{lat}(G)$, $\chi_{lat}(G-e)$, $\chi_{la}(G\vee K_1)$ and $\chi_{la}((G-e)\vee K_1)$ for many dense and sparse graphs $G$ in Section~\ref{sec-chilat}. We also conjecture that each graph $G$ of order at least 3 has $\chi_{lat}(G)\le \chi_{la}(G)$.

\section{Preliminary results}\label{sec-pre}

By definition, $\chi_{lat}(O_n) = \chi_{lat}(K_n) = n$. Moreover, if $G$ is a graph with $n$ isolated vertices, then $\chi_{lat}(G) \ge n$. In what follows, we only consider nonempty graphs.  

\begin{theorem}\label{thm-LAT} Every graph $G$ is local antimagic total. \end{theorem}

\begin{proof} For $G \vee K_1$ and $p\ge 1$, we let $G$ be of order $p$ and size $q$ such that $V(G)=\{v_i\,|\,1\le i\le p\}$, $V(K_1)=\{v\}$ and $E(G \vee K_1) = E(G)\cup \{vv_i \,|\, i = 1, 2, \ldots, p\}$. It is obvious that each graph $G$ of order $p\le 3$ are local antimagic total. We now assume $G$ is of order $p\ge4$. In~\cite{Haslegrave}, the author proved that every graph without isolated edges (by definition, necessarily without isolated vertices) admits a local antimagic labeling. Thus, $G \vee K_1$ is local antimagic. Let $g$ be a local antimagic labeling of $G \vee K_1$. Define a total labeling $f : V(G) \cup E(G) \to [1, p+q]$ of $G$ such that $f(e) = g(e)$ for each edge $e\in E(G)$ and $f(v_i) = g(vv_i)$. Clearly, $w(v_i) = g^+(v_i)$. Thus, $w(v_i)\neq w(v_j)$ if $v_iv_j\in E(G)$ and $f$ is a local antimagic total labeling of $G$.
\end{proof}

\nt In~\cite{Slamin+NMDK}, the authors extended the concept of local antimagic total labelings to local super antimagic total labelings (and chromatic number, denoted $\chi_{lsat}$) that assign only integers in $[1,|V(G)|]$ to the vertices of $G$. Since a local super antimagic total labeling is also a local antimagic labeling, but the converse may not hold, we then have $\chi(G)\le \chi_{lat}(G)\le \chi_{lsat}(G)$. By an argument similar to that of Theorem~\ref{thm-LAT}, we get a sufficient condition for a graph $G$ to admit a local super antimagic total labeling.

\begin{theorem} If $G\vee K_1$ admits a local antimagic labeling that assigns only integer(s) in $[1,|V(G)|]$ to the edge(s) not belong to $G$, then $G$ is local super antimagic total. \end{theorem}  

\nt The following three theorems follows directly from definition and Theorems 1, 7, 8 in~\cite{Slamin+NMDK} respectively.

\begin{theorem} If $T$ is a tree on $n \ge 2$ vertices with $k$ pendants, then $\chi_{lat}(T) \le n - k + 1$.  \end{theorem} 

\nt For graphs $G$ and $H$, the corona product of $G$ and $H$, denoted $G\odot H$, is obtained from $G$ and $|V(G)|\ge 2$ copies of $H$ by joining the $i$-th vertex of $G$ to every vertex of the $i$-th copy of $H$. 

\begin{theorem}\label{thm-HodotOm} If $H$ is a regular local super antimagic total graph of order $n \ge 2$ and $O_m$ is a null graph of order $m \ge 2$, then $\chi_{lat}(H \odot O_m) \le \chi_{lsat}(H) + 1$, where $(m, n) \not= (odd, even)$. \end{theorem}


\begin{theorem}\label{thm-HvK} Suppose $H$ is a regular graph of order $n \ge 2$ and $G$ is obtained from $H\vee K_1$ by joining each vertex of $H$ to $m$ pendant vertices for $m \ge 1$. If $G$ is local super antimagic total, then
\begin{enumerate}[(i)]
  \item $\chi_{lat}(G) \le \chi_{lsat}(H \vee K_1) + 1$, for $m \ge 2$ and $(m, n) \not= (odd, even)$,
  \item $\chi_{lat}(G) \le \chi_{lsat}(H) + 2$, for $(m, n) = (odd, even)$,
  \item $\chi_{lat}(G) \le \chi_{lsat}(H) + 2$, for $m = 1$.
\end{enumerate} 
\end{theorem}




\nt The next theorem shows that $\chi_{lat}(G)$ can be arbitrarily large for a graph $G$ with small $\chi(G)$.

\begin{theorem} If $G = K_2 + O_n, n\ge 1$, then $$\chi_{lat}(G)= \begin{cases} 2 &\mbox{ for } n=1,2, \\ n &\mbox{ otherwise.} \end{cases}$$
\end{theorem}

\begin{proof} Let $V(G) = \{u_1,u_2, v_i\,|\, 1\le i\le n\}$. For $n\ge 1$, define $f(u_i)=i$, $f(u_1u_2)=3$ and $f(v_i) = i+3$. We now have $w(u_1) = 4, w(u_2) = 5$ and $w(v_i) = i+3$. Thus, $\chi_{lat}(G) \le 2$ for $n=1,2$, and $\chi_{lat}(G)\le n$ for $n\ge 3$. By definition, $\chi_{lat}(G)\ge\chi(G)=  2$ and since all the isolated vertices must have distinct weights, this implies that $\chi_{lat}(G)\ge n$. So, the theorem holds. \end{proof}

\begin{theorem}\label{thm-K1VG} Let $G \vee K_1$ be as defined under Theorem~\ref{thm-LAT}. 
\begin{enumerate}[(a)]
  \item $\chi_{lat}(G) \le \chi_{la}(G \vee K_1) - 1$ and the equality holds if $\chi(G) = \chi_{la}(G \vee K_1)-1$. 
  \item Suppose $\chi_{lat}(G) = \chi(G \vee K_1) -1$ with a corresponding local antimagic total labeling $f$ . If $\sum_{i=1}^{p} f(v_i)\ne w(v_j), 1\le j\le p$, then $\chi_{la}(G \vee K_1)=\chi(G \vee K_1)$. \end{enumerate} \end{theorem}

\begin{proof}\hspace*{\fill}{}
 \begin{enumerate}[(a)]
   \item  Suppose $\chi_{la}(G \vee K_1)=c$. From the proof of Theorem~\ref{thm-LAT}, we know that every local antimagic labeling of $G \vee K_1$ that induces $c$ distinct vertex labels corresponds to a local antimagic total labeling of $G$ that induces $c-1$ distinct vertex weights. Thus, $\chi_{lat}(G)\le c-1$. By definition, $\chi_{lat}(G)\ge \chi(G)$. Thus, $\chi_{lat}(G) = c - 1$ if $\chi(G) = c-1$. 
   \item  Suppose $\chi(G \vee K_1) = c$ so that $\chi_{lat}(G)=c-1$. Let $g : E(G \vee K_1) \to [1,p+q]$ such that $g(e) = f(e)$ if $e\in E(G)$, and $g(vv_i) = f(v_i)$ for each $v_i\in V(G)$. Clearly, $w(v) = \sum_{i=1}^{p} f(v_i)$, and $w(v_i) = f^+(v_i)\ne f^+(v_j) = w(v_j) $ if $v_iv_j\in E(G)$. Since $w(v) \ne w(v_j)$, $g$ is a local antimagic labeling that induces $c$ distinct vertex weights and $\chi_{la}(G \vee K_1)\le c$. Since $\chi_{la}(G \vee K_1) \ge \chi(G \vee K_1) = c$, we have $\chi_{la}(G \vee K_1)=c$. \end{enumerate}
\end{proof}

\nt By an argument similar to that for Theorem~\ref{thm-K1VG}, we have the following theorem.

\begin{theorem}\label{cor-lsat} Let $G\vee K_1$ be as defined under Theorem~\ref{thm-LAT}. Suppose $G\vee K_1$ admits a local antimagic labeling that assigns only integers in $[1,p]$ to edges not belong to $G$.
\begin{enumerate}[(a)]
  \item $\chi_{lsat}(G)\le \chi_{la}(G \vee K_1) - 1$ and the equality holds if $\chi(G) = \chi_{la}(G \vee K_1)-1$. 
  \item Suppose $\chi_{lsat}(G) = \chi(G \vee K_1) -1$ with a corresponding local super antimagic total labeling $f$ . If ${p+1\choose2}\ne w(v_j), 1\le j\le p$, then $\chi_{la}(G \vee K_1)=\chi(G \vee K_1)$. \end{enumerate} \end{theorem}


\begin{corollary}\label{cor-K1veeG} For any graph $G$ of order at least 2,
\begin{enumerate}[(1)]
  \item $\chi(G) \le \chi_{lat}(G) \le \chi_{la}(G \vee K_1)-1$;
  \item if $G\vee K_1$ admits a local antimagic labeling that assigns only integers in $[1,p]$ to edges not belong to $G$, then $\chi(G) \le \chi_{lsat}(G) \le \chi_{la}(G \vee K_1)-1$.  
\end{enumerate}  
\end{corollary}

\nt For $1\le i\le p$ and $1\le j\le m$,  let $G$ be a graph of order $p\ge 2$ and size $q\ge 1$ with $V(G \vee mK_1) = \{v_i, u_j\}$ and $E(G\vee mK_1)=E(G)\cup \{u_iv_j\}$. 

\begin{theorem}\label{thm-GvO2} Suppose $g$ is a local antimagic labeling of $G\vee 2K_1$ that induces a minimum number of vertex labels and $2p+q+1+g^+(u_1) \ne g^+(v_i), 1\le i\le p$, then $$\chi_{lat}(G \vee K_1) \le \begin{cases} \chi_{la}(G\vee 2K_1) & \mbox{ if } g^+(u_1) = g^+(u_2); \\ \chi_{la}(G\vee 2K_1) - 1 & \mbox{ if } g^+(u_1) \ne g^+(u_2).   \end{cases}$$
Moreover, equality holds if $\chi(G \vee K_1)$ equals the respective upper bound.
\end{theorem}

\begin{proof} For $1\le i < j\le p$, define a bijection $f:V(G \vee K_1)\cup E(G \vee K_1) \to [1,2p+q+1]$ such that for $1\le i < j \le p$, $f(v_iv_j) = g(v_iv_j)$ if $v_iv_j\in E(G)$, $f(u_1v_i) = g(u_1v_i)$, $f(v_i) = g(u_2v_i)$, and $f(u_1)=2p+q+1$. Clearly, $w(v_i) = g^+(v_i)$ and $w(u_1) = 2p+q+1+g^+(u_1)$. Since $2p+q+1+g^+(u_1) \ne g^+(v_i)$ for $1\le i\le p$, $f$ is a local antimagic total labeling of $G \vee K_1$ that induces $\chi_{la}(G\vee 2K_1)$ distinct vertex weights if $g^+(u_1) = g^+(u_2)$, and induces $\chi_{la}(G\vee 2K_1) - 1$ distinct vertex weights if $g^+(u_1) \ne g^+(u_2)$. Moreover, equality holds if $\chi(G \vee K_1)$ equals the respective upper bound.
\end{proof}

\nt 
The following lemmas are analogous to Lemmas 2.2 -- 2.5 in~\cite{LSN1}.

\begin{lemma}\label{lem-regular} Suppose $G$ is a $d$-regular graph of order $p$ and size $q$ with an edge $e$. If $f$ is a local antimagic total labeling of $G$, then $g=p+q+1-f$ is also a local antimagic total labeling of $G$ with $w(g) = w(f)$. Moreover, suppose $w(f) = \chi_{lat}(G)$ and $f(e)=1$ or $f(e) = p+q$. If $\chi(G-e)=\chi_{lat}(G)$, then $\chi_{lat}(G - e) = \chi_{lat}(G)$. Otherwise, $\chi_{lat}(G-e)\le \chi_{lat}(G)$. \end{lemma}

\begin{proof} Let $x,y\in V(G)$. Here, $w_g(x) = (d+1)(p+q+1) - w_f(x)$ and $w_g(y) = (d+1)(p+q+1) - w_f(y)$. Therefore, $w_f(x) = w_f(y)$ if and only if $w_g(x) = w_g(y)$. Thus, $g$ is also a local antimagic total labeling of $G$ with $w(g) = w(f)$.

\ms\nt If $f(e) = p+q$, then we may consider $g=p+q+1-f$. So without loss of generality, we may assume that $f(e)=1$. Define $h:V(G-e)\cup E(G-e)\to [1,p+q-1]$ such that $h(x) = f(x) - 1$ and $h(xy) = f(xy)-1$ for $xy\ne e$. So, $w_h(x) = w_f(x) - d-1$ for each vertex $x$ of $G - e$. Therefore, $w_f(x) = w_f(y)$ if and only if $w_h(x) = w_h(y)$. Thus, $h$ is also a local antimagic total labeling of $G$ with $w(h) = w(f)$. Consequently, $\chi(G-e)\le \chi_{lat}(G-e) \le \chi_{lat}(G)$.  The theorem holds.
\end{proof}

\nt Note that if $G$ is a regular edge-transitive graph, then $\chi_{lat}(G - e) \le \chi_{lat}(G)$.

\begin{lemma}\label{lem-non-regular} Suppose $G$ is a graph of order $p$ and size $q$ and $f$ is a local antimagic total labeling of $G$. For any $x,y\in V(G)$, if
\begin{enumerate}[(i)]
  \item $w_f(x) = w_f(y)$ implies that $deg(x)=deg(y)$, and
  \item $w_f(x)\ne w_f(y)$ implies that $(p+q+1)(deg(x) - deg(y)) \ne w_f(x) - w_f(y)$,   
\end{enumerate}
then $g=p+q+1-f$ is also a local antimagic total labeling of $G$ with $w(g) = w(f)$. 
\end{lemma}

\begin{proof} For any $x,y\in V(G)$, we have $w_g(x) = (deg(x)+1)(p+q+1)-w_f(x)$ and $w_g(y) = (deg(y)+1)(p+q+1)-w_f(y)$. If $w_f(x) = w_f(y)$, then condition (i) implies that $w_g(x) = w_g(y)$. If $w_f(x) \ne w_f(y)$, then condition (ii) implies that $w_g(x) \ne w_g(y)$. Thus, $g$ is also a local antimagic total labeling of $G$ with $w(g)=w(f)$.
\end{proof}

\nt For $t\ge 2$, consider the following conditions for a graph $G$. 

\begin{enumerate}[(i)]
  \item $\chi_{lat}(G) = t$ and $f$ is a local antimagic total labeling of $G$ that induces a $t$-independent partition $\bigcup^t_{i=1} V_i$ of $V(G)$.
  \item For each $x\in V_k$, $1\le k\le t$, $deg(x) = d_k$ satisfying $w_f(x) - d_a \ne w_f(y) - d_b$, where $x\in V_a$ and $y\in V_b$ for $1\le a < b\le t$.
  \item There exist two non-adjacent vertices $u,v$ with $u\in V_i, v\in V_j$ for some $1\le i\ne j\le t$ such that 
  \begin{enumerate}[(a)]
    \item $|V_i| = |V_j| = 1$ and $deg(x) = d_k$ for $x\in V_k, 1\le k\le t$; or
    \item $|V_i| = 1, |V_j|\ge 2$ and $deg(x) = d_k$ for $x\in V_k, 1\le k\le t$ except that $deg(v) = d_j-1$; or
    \item $|v_i|, |V_j|\ge 2$ and $deg(x) = d_k$ for $x\in V_k$, $1\le k\le t$ except that $deg(u) = d_i-1$, $deg(v) = d_j - 1$,
  \end{enumerate} 
  each satisfying $w_f(x) + d_a \ne w_f(y) + d_b$, where $x\in V_a$ and $y\in V_b$ for $1\le a\ne b \le t$.
\end{enumerate}

\begin{lemma}\label{lem-lat-G-e} Let $H$ be obtained from $G$ with an edge $e$ deleted. If $G$ satisfies conditions (i) and $(ii)$ and $f(e) = 1$, then $\chi(H)\le \chi_{lat}(H)\le t$.   \end{lemma}

\begin{proof} By definition, we have the lower bound. Define $g: E(H) \to [1,|E(H)|]$ such that $g(e') = f(e') - 1$ for each $e'\in E(H)$. Observe that $g$ is a bijection with $w_g(x) = w_f(x) - d_k - 1$ for each $x\in V_k, 1\le k\le t$. Thus, $w_g(x) = w_g(y)$ if and only if $x,y\in V_k, 1\le k\le t$. Therefore, $g$ is a local antimagic total labeling of $H$ with $w(g) = w(f)$. Thus, $\chi_{lat}(H)\le t$.
\end{proof}

\begin{lemma}\label{lem-lat-G+e} Suppose $uv\in E(G)$. Let $H$ be obtained from $G$ with an edge $uv$ added. If $G$ satisfies conditions (i) and (iii), then $\chi(H)\le \chi_{lat}(H) \le t$. \end{lemma}

\begin{proof} By definition, we have the lower bound. Define $g : E(H) \to [1,|E(H)|]$ such that $g(uv)=1$ and $g(e) = f(e) + 1$ for $e \in E(G)$. Observe that $g$ is a bijection with $w_g(x) = w_f(x) + d_k + 1$ for each $x \in V_k, 1\le k\le t$. Thus, $w_g(x) = w_g(y)$ if and only if $x,y\in V_k, 1\le k\le t$. Therefore, $g$ is a local antimagic total labeling of $H$ with $w(g) = w(f)$. Thus, $\chi_{lat}(H) \le t$. \end{proof}

\section{The exact local antimagic total chromatic number}\label{sec-chilat}

\nt 
We first consider some dense graphs. An amalgamation of two graphs $G_1$ and $G_2$ over a fixed graph $H$ is the simple graph obtained by identifying the vertices of two induced subgraphs isomorphic to $H$, one of $G_1$ and the other of $G_2$. Suppose $G$ is a graph with a $K_r$, $r\ge 1$, subgraph. Let $A(mG,K_r)$ be the amalgamation of $m\ge 2$ copies of $G$ of order at least $r+1$ along $K_r$. When $r=1$, the graph is also known as one-point union of graphs. For $m\ge2$ and $n > r\ge 1$, let $V(A(mK_n,K_r))=\{v_{i,j}\,|\,1\le i\le m, 1\le j\le n\}$ and $E(A(mK_n,K_r)) = \{v_{i,j}v_{i,k}\,|\, 1\le i\le m, 1\le j < k\le n\}$ where $v_{1,j}=\cdots = v_{m,j}$ for each $1\le j\le r$. Note that $A(mK_2,K_1)\cong K_{1,m}$ and $A(mK_3, K_1)$ is the friendship graph $f_m, m\ge 2$. 


\begin{theorem}\label{thm-A(mKnKr)} For $m\ge 3, n\ge r+1\ge 3$, if $3\le m\le 4$ and $2\le r\le 3$, or if $(m-1)(r-1) < {r+1\choose2}$, then $\chi_{lat}(A(mK_n,K_r)) = n$. \end{theorem}

\begin{proof} Let $V(A(mK_n,K_r)) = \{v_{i,j}\,|\,1\le i\le m, 1\le j\le n\}$ and $E(A(mK_n,K_r)) = \{v_{i,j}v_{i,k}\,|\,1\le i\le m, 1\le j<k\le n\}$ such that $\{v_{i,j}\,|\,1\le j\le r\}$ is the set of vertices common to each $K_n$. Thus, $v_{1,j}=v_{2,j}=\cdots=v_{m,j}$ for $1\le j\le r$. Define $f:V(A(mK_n,K_r))\cup E(A(mK_n,K_r)) \to [1,m{n+1\choose2}-(m-1){r+1\choose2}]$ as follows.
\begin{enumerate}[(i)]
  \item $f(v_{i,j}v_{i,k})=m\Big[{k-1\choose2}-{r\choose2}+j-1\Big]+i$ for $1\le i\le m, 1\le j < k\le n, r+1 \le k\le n$;
  \item $f(v_{i,j}) = m({n\choose2}-{r\choose2})+(m+1-i)(n-1)+2-j$ for $1\le i\le m, r+1\le j\le n$. 
\end{enumerate}


\ms\nt Let $t=m({n\choose2}-{r\choose2})$. We now consider the labeling of the remaining $r$ vertices and ${r\choose2}$ edges that belong to the $K_r$. 

\nt {\bf Case (1).} Suppose $m=4$. If $r=2$, we are left with $t+n-1,t+2n-2,t+3n-3$ to label the vertices and edge of the unlabeled $K_2$ subgraph. Let $f(v_{i,2})=t+n-1$, $f(v_{i,1})=t+2n-2$ and $f(v_{i,1}v_{i,2})=t+3n-3$. We now have the followings.

\begin{eqnarray*}
w(v_{i,1})&=&\sum^4_{i=1}\sum^n_{k=3}\Bigg\{4\Big[{k-1\choose2}-1\Big]+i\Bigg\}+2t+5n-5\\
&=&16\Bigg[{n\choose3}-n+2\Bigg]+10(n-2)+8\Bigg[{n\choose2}-1\Bigg]+5n-5\\
&=&\frac{8n^3}{3}-4n^2+\frac{n}{3}-1
\end{eqnarray*}

and 
\begin{eqnarray*}
w(v_{i,2})&=&\sum^4_{i=1}\sum^n_{k=3}\Bigg\{4{k-1\choose2}+i\Bigg\}+2t+4n-4\\
&=&16{n\choose3}+10(n-2)+8\Bigg[{n\choose2}-1\Bigg]+4n-4\\
&=&\frac{8n^3}{3}-4n^2+\frac{46n}{3}-32.
\end{eqnarray*}

\nt For $3\le j\le n$, 
\begin{eqnarray*}
w(v_{i,j})&=& f(v_{i,j}) + \sum^{j-1}_{k=1}f(v_{i,j}v_{i,k}) + \sum^n_{k=j+1}f(v_{i,j}v_{i,k})\\
&=& 4{n\choose2}-4+(5-i)(n-1)+2-j + \sum^{j-1}_{k=1}\Bigg[4\Bigg({j-1\choose2}+k-2\Bigg)+i\Bigg]+\\
& & \sum^n_{k=j+1}\Bigg[4\Bigg({k-1\choose2}+j-2\Bigg)+i\Bigg]\\
&=& 4{n\choose2}-4 +5(n-1)+2-j+ 4(j-1)\Bigg[{j-1\choose2}-2\Bigg]+4{j\choose2}+\\
& & 4(n-j)(j-2)+4\Bigg({n\choose3}-{j\choose3}\Bigg)\\
&=& \frac{2n^3}{3}+4nj-\frac{11n}{3}+\frac{4j^3}{3}-8j^2+\frac{17j}{3}-3.
\end{eqnarray*}
which is an increasing function. Moreover,

\begin{eqnarray*}
w(v_{i,1}) - w(v_{i,n})&=&\frac{8n^3}{3}-4n^2+\frac{n}{3}-1-
\Bigg(\frac{2n^3}{3}+4n^2-\frac{11n}{3}+\frac{4n^3}{3}-8n^2+\frac{17n}{3}-3\Bigg)\\
&=& \frac{2n^3}{3}-\frac{5n}{3}+2>0.
\end{eqnarray*}

\nt Since $w(v_{i,3}) < \cdots < w(v_{i,n}) <w(v_{i,1}) < w(v_{i,2})$, we have $f$ is a local antimagic total labeling that induces $n$ distinct vertex weights so that $\chi_{la}(A(mK_n,K_r))\le n$.

\ms\nt If $r=3$, we are left with $t+n-2,t+n-1,t+2n-3,t+2n-2,t+3n-4,t+3n-3$ to label the vertices and edges of the unlabeled $K_3$ subgraph. Let $f(v_{i,3})=t+n-2$, $f(v_{i,2})=t+n-1$, $f(v_{i,2}v_{i,3})=t+2n-3$, $f(v_{i,1})=t+2n-2$, $f(v_{i,1}v_{1,2})=t+3n-4$ and $f(v_{i,1}v_{1,3})=t+3n-3$. We now have the followings.

\begin{eqnarray*}
w(v_{i,1})&=& \sum^4_{i=1}\sum^n_{k=4}\Bigg\{4\Big[{k-1\choose2}-3\Big]+i\Bigg\} + 3t+8n-9 \\
&=&\frac{8n^3}{3}-2n^2-\frac{92n}{3}+53,
\end{eqnarray*}

\begin{eqnarray*}
w(v_{i,2})&=& \sum^4_{i=1}\sum^n_{k=4}\Bigg\{4\Big[{k-1\choose2}-2\Big]+i\Bigg\} + 3t+6n-8\\
&=&\frac{8n^3}{3}-2n^2-\frac{50n}{3}+6
\end{eqnarray*}

and 
\begin{eqnarray*}
w(v_{i,3})&=& \sum^4_{i=1}\sum^n_{k=4}\Bigg\{4\Big[{k-1\choose2}-1\Big]+i\Bigg\} + 3t+6n-8\\
&=&\frac{8n^3}{3}-2n^2-\frac{2n}{3}-42.
\end{eqnarray*}

\nt For $4\le j\le n$, 
\begin{eqnarray*}
w(v_{i,j})&=& f(v_{i,j}) + \sum^{j-1}_{k=1}f(v_{i,j}v_{i,k}) + \sum^n_{k=j+1}f(v_{i,j}v_{i,k})\\
&=& \frac{2n^3}{3}+4nj-\frac{35n}{3}+\frac{4j^3}{3}-8j^2+\frac{17j}{3}-3
\end{eqnarray*}
which is an increasing function. Moreover,

 \begin{eqnarray*}
w(v_{i,1})-w(v_{i,n})&=&\frac{8n^3}{3}-2n^2-\frac{92n}{3}+53-\Big(\frac{2n^3}{3}+4n^2-\frac{35n}{3}+\frac{4n^3}{3}-8n^2+\frac{17n}{3}-3\Big)\\
&=&\frac{2n^3}{3}+2n^2-\frac{74n}{3}+56>0.
\end{eqnarray*}

\nt Since $w(v_{i,4}) < \cdots < w(v_{i,n}) < w(v_{i,1}) < w(v_{i,2}) < w(v_{i,3})$, we have $f$ is a local antimagic total labeling that induces $n$ distinct vertex weights so that $\chi_{la}(A(mK_n,K_r))\le n$.

\ms\nt {\bf Case (2).} Suppose $m=3$. The proof is similar to Case (1) and is thus omitted. We note that when $r=2$, we can get $w(v_{i,1})=\frac{3n^2(n-1)}{2}+2n-5$, $w(v_{i,2})=\frac{3n^2(n-1)}{2}+10n-22$, and for $3\le j\le $, $w(v_{1,j})=\frac{n^3}{2} +3nj-\frac{5n}{2}+j^3-6j^2+4j-2$, which is an increasing function with $w(v_{i,1}) - w(v_{i,n})=\frac{2n^3}{3}-\frac{5n}{3}+2>0$. For $r=3$, we can get $w(v_{i,1})=\frac{3n^3}{2}-\frac{29n}{2}+18$, $w(v_{i,2})=\frac{3n^3}{2}-\frac{15n}{2}-8$, $w(v_{i,3})=\frac{3n^3}{2}+\frac{3n}{2}-35$, and for $4\le j\le n$, $w(v_{i,j})=\frac{n^3}{2}+3nj-\frac{17n}{2}+j^3-6j^2+4j-2$, which is an  increasing function with $w(v_{i,1})-w(v_{i,n})=3n^2-10n+20>0$.

\ms\nt {\bf Case (3).} Suppose $(m-1)(r-1) < {r+1\choose2}$. We may assume $(m,r)\ne (4,2),(4,3),(3,2),(3,3)$. Note that the ${r+1\choose2}$ unused integers in $[1,m{n+1\choose2}-(m-1){r+1\choose2}]$ are as in 

\begin{enumerate}[(i)]
  \item $a_i=\{t+in-r+2-i, t+in-r+3-i, t+in-r+4-i, \ldots, t+in-i\}$ that has $r-1$ integers for $1\le i\le m$, and 
  \item $a_{m+1}=\{t+mn-m+1, t+mn-m+2, \ldots, t+mn+{r\choose2}-r(m-1)\}$ that has ${r\choose2}+m-r(m-1)$ integers. 
\end{enumerate}

\nt Let us list the above integers in increasing order to get \\
\nt $A(1) = t+n-r+1, A(2)=t+n-r+2, A(3)=t+n-r+3, \ldots, A(r-1)=t+n-1,$ \\$A(r)=t+2n-r, A(r+1)=t+2n-r+1, A(r+2)=t+2n-r+2, \ldots, A(2r-2)=t+2n-2,$ \\ $A(2r-1)=t+3n-r-1, A(2r)=t+3n-r, A(2r+1)=t+3n-r+1, \ldots, A(3r-3) = t+3n-3,$ \\ $\vdots$ \\ $A((m-1)(r-1)+1)=t+mn-r-m+2, A((m-1)(r-1)+2)=t+mn-r-m+3, A((m-1)(r-1)+3)=t+mn-r-m+4, \ldots, A(m(r-1))=t+mn-m,$ and \\ $A(m(r-1)+1)=t+mn-m+1, A(m(r-1)+2)=t+mn-m+2, \ldots, A({r+1\choose2})=t+mn+{r\choose2}-r(m-1)$.\\

\nt Let $M_{r\times r}$ be an $r\times r$ symmetric matrix that corresponds to the labeling of the $K_r$ subgraph such that the integers $A(1)$ to $A(r)$ are assigned to the main diagonal from bottom right to top left, and the remaining integers from $A(r+1)$ onwards are assigned to $M_{j,j+1}$ to $M_{j,r}$ for $j=1$ to $j=r$ consecutively. Thus, the row sums of $M_{r\times r}$ that corresponds to the total labeling of the $K_r$-subgraph are:\\
\nt $S_1 = {\bf A(r)} + A(r+1) + A(r+2) + A(r+3) + A(r+4) + A(r+5) + A(r+6) + A(r+7) + \ldots + A(2r-1)$;\\
\nt $S_2 = A(r+1) + {\bf A(r-1)} + A(2r) + A(2r+1) + A(2r+2) + A(2r+3) + A(2r+4) + A(2r+5) + \ldots + A(3r-3)$;\\
\nt $S_3 = A(r+2) + A(2r) + {\bf A(r-2)} + A(3r-2) + A(3r-1) + A(3r) + A(3r+1) + A(3r+2) + \ldots + A(4r-6)$;\\
\nt $S_4 = A(r+3) + A(2r+1) + A(3r-2) + {\bf A(r-3)}+ A(4r-5) + A(4r-4) + A(4r-3) + A(4r-2) +\ldots + A(5r-10)$;\\
\nt $S_5 = A(r+4) + A(2r+2) + A(3r-1) + A(4r-5) + {\bf A(r-4)} + A(5r-9) + A(5r-8) + A(5r-7) +  \ldots + A(6r-15)$;\\
\nt $S_6 = A(r+5) + A(2r+3) + A(3r) + A(4t-4) + A(5r-9) + {\bf A(r-5)} + A(6r-14) + A(6r-13) + \ldots + A(7r-21)$;\\
$\vdots$\\
\nt $S_r = A(1) + \sum^{r-1}_{x=1}A((x+1)r-{x+1\choose2})$.\\
\nt In general, for $2\le x\le r-1$, we have \\ $S_x = A(r+x-1) + A(2r+x-3) + A(3r+x-5) + \ldots + A((x-1)r-{x-1\choose2}+1) + {\bf A(r-x+1)} + A(xr-{x\choose2}+1) + A(xr-{x\choose2}+2) + \ldots + A((x+1)r-{x+1\choose2})$.
\ms\nt It is now obvious that $S_1 < S_2 < \cdots < S_r$. Moreover, 
$$S_1 = t+3n-r-1+\sum^{r-1}_{y=1} (t+2n-r-1+y) = r(t+2n-1)-{r+1\choose2}+n.$$ Note that for given $m,n,r$,
\begin{eqnarray*}
\sum^{m}_{i=1}\sum^{n}_{k=r+1} f(v_{i,j}v_{i,k}) &=& \sum^{m}_{i=1}\sum^{n}_{k=r+1} \Bigg(m\Big[{k-1\choose2}-{r\choose2}+j-1\Big]+i\Bigg)\\
&=&m^2\Big({n\choose3}-{r\choose3}\Big) + (n-r){m+1\choose2}+m^2(n-r)\Big[j-1-{r\choose2}\Big]
\end{eqnarray*} 
is an increasing function from $j=1$ to $j=r$.  Thus, we can now conclude that $w(v_{i,1}) < w(v_{i,2}) <\cdots < w(v_{i,r})$ and 
\begin{eqnarray*} w(v_{i,1}) &=& S_1 + m^2\Big({n\choose3}-{r\choose3}\Big) + (n-r){m+1\choose2}-m^2(n-r){r\choose2}\\
&=&r(t+2n-1)-{r+1\choose2}+n + m^2\Big({n\choose3}-{r\choose3}\Big) + (n-r){m+1\choose2}-m^2(n-r){r\choose2}
\end{eqnarray*}

\ms\nt Also note that for $r+1\le j\le n$, we have
\begin{eqnarray*}
w(v_{i,j})&=&f(v_{i,j})+\sum^{j-1}_{k=1}f(v_{i,j}v_{i,k}) + \sum^n_{k=j+1}f(v_{i,j}v_{i,k})\\
&=& t+(m+1-i)(n-1)+2-j + \sum^{j-1}_{k=1}\Bigg[ m\Big[{j-1\choose2}-{r\choose2}+k-1\Big]+i\Bigg]\\
& & \sum^n_{k=j+1} \Bigg[ m\Big[{k-1\choose2}-{r\choose2}+j-1\Big]+i\Bigg]\\
&=& t+(m+1)(n-1)+2-j + m(j-1)\Big[{j-1\choose2}-{r\choose2}-1\Big]+m{j\choose2}+ \\
& & m\Big[{n\choose3}-{j\choose3}\Big]+m(n-j)\Big[j-1-{r\choose2}\Big].\\
\end{eqnarray*}

\nt It is routine to check that the terms involving the variable $j$ form an increasing function so that $w(v_{i,r+1}) < w(v_{i,r+2}) < \cdots < w(v_{i,n})$. Since $n\ge r+1$, we get

\begin{eqnarray*}
w(v_{i,1}) - w(v_{i,n})&=&r(t+2n-1)-{r+1\choose2}+n + m^2\Big({n\choose3}-{r\choose3}\Big) + (n-r){m+1\choose2}-\\
&& m^2(n-r){r\choose2}-\Bigg(t+(m+1)(n-1)+2-n + m(n-1)\Big[{n-1\choose2}-{r\choose2}-1\Big]+\\
&& m{n\choose2}\Bigg)\\
&=&\frac{1}{6}\Big[2m^2r^3-3mr^3-3m^2nr^2+3mnr^2+3mr^2-3r^2+3mn^2r+3m^2nr-\\
&& 6mnr+12nr-5m^2r-3mr-9r+m^2n^3-3mn^3-3m^2n^2+6mn^2+\\
&& 5m^2n-6mn+6n+6m-6\Big]\\
&>&0.
\end{eqnarray*}
Thus, $w(v_{i,r+1}) < w(v_{i,r+2}) <\cdots < w(v_{i,n}) < w(v_{i,1}) < \cdots < w(v_{i,r})$ so that $f$ is a local antimagic total labeling that induces $n$ distinct vertex weights. Therefore, $\chi_{lat}(A(mK_n,K_r))\le n$. Since $\chi_{la}(A(mK_n,K_r))\ge \chi(A(mK_n,K_r))=n$, the theorem holds.
\end{proof} 

\begin{rem} We note that when $m=2$, we may have $w(v_{i,1}) - w(v_{i,n})<0$. We have checked by computer computation that all the induced vertex weights are distinct for $3\le r+1\le n\le 15000$. \end{rem}

\begin{example} For Cases (1) and (2), we give the labeling matrix of $A(4K_7,K_3)$ denoted $A(i,7,3)$ for $1\le i\le 4$. \\

\nt\begin{minipage}[b]{9cm}
$\fontsize{9}{10}\selectfont\begin{array}{c||*{7}{c|}|c}
A(1,7,3) & v_{1,1} & v_{1,2} & v_{1,3} & v_{1,4} & v_{1,5} & v_{1,6} & v_{1,7} & \mbox{sum} \\\hline\hline
v_{1,1} & 84 & 89 & 90 & 1 & 13 & 29 & 49 & 355 \\\hline
v_{1,2} & 89 & 78 & 83 & 5 & 17 & 33 & 53 & 358 \\\hline
v_{1,3} & 90 & 83 & 77 & 9 & 21 & 37 & 57 & 374 \\\hline
v_{1,4} & 1 & 5 & 9 & 94 & 25 & 41 & 61 & 236 \\\hline
v_{1,5} & 13 & 17 & 21 & 25 & 93 & 45 &  65 & 279 \\\hline
v_{1,6} & 29 & 33 & 37 & 41 & 45 & 92 & 69 & 346 \\\hline
v_{1,7} & 49 & 53 & 57 & 61 & 65 & 69 & 91 & 445 \\\hline
\end{array}$\\
\end{minipage}
\begin{minipage}[b]{9cm}
$\fontsize{9}{10}\selectfont\begin{array}{c||*{7}{c|}|c}
A(2,7,3) & v_{2,1} & v_{2,2} & v_{2,3} & v_{2,4} & v_{2,5} & v_{2,6} & v_{2,7} & \mbox{sum}\\\hline\hline
v_{2,1} & - & - & - & 2 & 14 & 30 & 50 & 96 \\\hline
v_{2,2} & - & - & - & 6 & 18 & 34 & 54 & 112\\\hline
v_{2,3} & - & - & - & 10 & 22 & 38 & 58 & 128 \\\hline
v_{2,4} & 2 & 6 & 10 & 88 & 26 & 42 & 62 & 236 \\\hline
v_{2,5} & 14 & 18 & 22 & 26 & 87 & 46 & 66 & 279 \\\hline
v_{2,6} & 30 & 34 & 38 & 42 & 46 & 86 & 70 & 346 \\\hline
v_{2,7} & 50 & 54 & 58 & 62 & 66 & 70 & 85 & 445 \\\hline
\end{array}$\\
\end{minipage}\\\\
\nt\begin{minipage}[b]{9cm}
$\fontsize{9}{10}\selectfont\begin{array}{c||*{7}{c|}|c}
A(3,7,3) & v_{3,1} & v_{3,2} & v_{3,3} & v_{3,4} & v_{3,5} & v_{3,6} & v_{3,7} & \mbox{sum}\\\hline\hline
v_{3,1} & - & - & - & 3 & 15 & 31 & 51 & 100 \\\hline
v_{3,2} & - & - & - & 7 & 19 & 35 & 55 & 116 \\\hline
v_{3,3} & - & - & - & 11 & 23 & 39 & 59 & 132 \\\hline
v_{3,4} & 3 & 7 & 11 & 82 & 27 & 43 & 63 & 236 \\\hline
v_{3,5} & 15 & 19 & 23 & 27 & 81 & 47 & 67 & 279 \\\hline
v_{3,6} & 31 & 35 & 39 & 43 & 47 & 80 & 71 & 346 \\\hline
v_{3,7} & 51 & 55 & 59 & 63 & 67 & 71 & 79 & 445 \\\hline
\end{array}$\\
\end{minipage}
\begin{minipage}[b]{9cm}
$\fontsize{9}{10}\selectfont\begin{array}{c||*{7}{c|}|c}
A(4,7,3) &v_{4,1} & v_{4,2} & v_{4,3} & v_{4,4} & v_{4,5} & v_{4,6} & v_{4,7} & \mbox{sum}\\\hline\hline
v_{4,1} & - & - & - & 4 & 16 & 32 & 52 & 104 \\\hline
v_{4,2} & - & - & - & 8 & 20 & 36 & 56 & 120 \\\hline
v_{4,3} & - & - & - & 12 & 24 & 40 & 60 & 136 \\\hline
v_{4,4} & 4 & 8 & 12 & 76 & 28 & 44 & 64 & 236 \\\hline
v_{4,5} & 16 & 20 & 24 & 28 & 75 & 48 & 68 & 279 \\\hline
v_{4,6} & 32 & 36 & 40 & 44 & 48 & 74 & 72 & 346 \\\hline
v_{4,7} & 52 & 56 & 60 & 64 & 68 & 72 & 73 & 445 \\\hline
\end{array}$\\
\end{minipage}\\
\nt The vertex weights are $w(v_{i,1}) = 355+96+100+104=685$, $w(v_{i,2}) = 358+112+116+120=706$, $w(v_{i,3}) = 374+128+132+136=770$, $w(v_{i,4}) = 236$, $w(v_{i,5})=279$, $w(v_{i,6})=346$ and $w(v_{i,7})=445$ for $1\le i\le 4$.
\end{example}

\begin{example}  For $(m-1)(r-1)<{r+1\choose2}$, we give the labeling matrix of $A(2K_6,K_3)$ denoted $A(i,6,3)$ for $1\le i\le 2$. \\

\nt\begin{minipage}[b]{9cm}
$\fontsize{9}{10}\selectfont\begin{array}{c||*{6}{c|}|c}
A(1,6,3) & v_{1,1} & v_{1,2} & v_{1,3} & v_{1,4} & v_{1,5} & v_{1,6} & \mbox{sum} \\\hline\hline
v_{1,1} & 33 & 34 & 35 & 1 & 7 & 15 & 125 \\\hline
v_{1,2} & 34 & 29 & 36 & 3 & 9 & 17 & 128 \\\hline
v_{1,3} & 35 & 36 & 28 & 5 & 11 & 19 & 134 \\\hline
v_{1,4} & 1 & 3 & 5 & 32 & 13 & 21 & 75 \\\hline
v_{1,5} & 7 & 9 & 11 & 13 & 31 & 23 & 94 \\\hline
v_{1,6} & 15 & 17 & 19 & 21 & 23 & 30 & 125 \\\hline
\end{array}$\\
\end{minipage}
\begin{minipage}[b]{9cm}
$\fontsize{9}{10}\selectfont\begin{array}{c||*{6}{c|}|c}
A(2,6,3) & v_{2,1} & v_{2,2} & v_{2,3} & v_{2,4} & v_{2,5} & v_{2,6} & \mbox{sum}\\\hline\hline
v_{2,1} & - & - & - & 2 & 8 & 16 & 26 \\\hline
v_{2,2} & - & - & - & 4 & 10 & 18 & 32\\\hline
v_{2,3} & - & - & - & 6 & 12 & 20 & 38 \\\hline
v_{2,4} & 2 & 4 & 6 & 27 & 14 & 22 & 75 \\\hline
v_{2,5} & 8 & 10 & 12 & 14 & 26 & 24 & 94 \\\hline
v_{2,6} & 16 & 18 & 20 & 22 & 24 & 25 & 125 \\\hline
\end{array}$\\
\end{minipage}\\
\nt The vertex weights are $w(v_{i,1}) = 125+26=151$, $w(v_{i,2}) = 121+32=160$, $w(v_{i,3}) = 134+38=172$, $w(v_{i,4}) = 75$, $w(v_{i,5})=94$ and $w(v_{i,6})=125$ for $1\le i\le 2$.
\end{example}

\begin{theorem}\label{thm-A(mKnKr)vK1} For $m\ge2$ and $n\ge r+1\ge 2$,  if $3\le m\le 4$ and $2\le r\le 3$, or if $(m-1)(r-1)\le {r+1\choose2}$, then $\chi_{la}(A(mK_n,K_r)\vee K_1) = n+1$.
\end{theorem}

\begin{proof} Suppose $m=4$. Consider $r=2$ and the local antimagic total labeling $f$ in the proof of Theorem~\ref{thm-A(mKnKr)}. We have the sum of all vertex labels is
\begin{eqnarray*}
\sum^2_{j=1} f(v_{1,j}) + \sum^4_{i=1}\sum^n_{j=3} f(v_{i,j})&=&(t+2n-2)+(t+n-1)+\\
& &\sum^4_{i=1}\sum^n_{j=3}\Big(4{n\choose2}+(5-i)(n-1)-2-j\Big)\\
&=& 2t+16(n-2){n\choose2}-4{n+1\choose2}+10n^2-27n+21\\
&>&w(v_{i,2}).
\end{eqnarray*}
\nt Since $\chi_{lat}(A(4K_n,K_2))=n=\chi(A(4K_n,K_2)\vee K_1)-1$, by Theorem~\ref{thm-K1VG}, $\chi_{la}(A(4K_n,K_2)\vee K_1)=n+1$.

\ms\nt For all other cases when $m=3,4$, one can check that $f$ satisfies the conditions of Theorem~\ref{thm-K1VG}. The details are omitted.

\ms\nt Suppose $(m-1)(r-1) < {r+1\choose2}$. It is routine to check that the sum of all vertex labels is $S_v=m^2(n-r)\Big[{n\choose2}-{r\choose2}\Big] + (n-r)(n-1){m+1\choose2}+2m(n-r)-\frac{1}{2}m(n-r)(n+r+1)+mr\Big[{n\choose2}-{r\choose2}\Big]+n(r+1)-{r+1\choose2}$ while the smallest vertex labels is $w(v_{i,r+1}) = m\Big[{n\choose2}-{r\choose2}\Big]+(m+1)(n-1)-r-mr+m{r+1\choose2}+m\Big[{n\choose3}-{r+1\choose3}\Big]+m(n-r-1)\Big[r-{r\choose2}\Big]$. Thus, $S_v - w(v_{i,r+1}) = \frac{1}{2}m^2r^3-\frac{5}{6}mr^3-\frac{1}{2}m^2nr^2+\frac{1}{2}mnr^2-m^2r^2+2mr^2-\frac{r^2}{2}-\frac{1}{2}m^2n^2r+\frac{1}{2}mn^2r+m^2nr-\frac{5}{2}mnr+nr+\frac{1}{2}m^2r+\frac{1}{3}mr+\frac{r}{2}+\frac{1}{2}m^2n^3-\frac{1}{6}mn^3-\frac{1}{2}m^2n+\frac{mn}{6}+m+1>0$. Thus, the theorem holds. 
\end{proof}


\nt Let $e$ be an edge of $A(mK_n,K_r)$ joining a vertex of the $K_r$ subgraph and another vertex not in $K_r$.

\begin{theorem}\label{thm-A(mKnKr-e)} For $m\ge2$ and $n\ge r+1\ge 2$,  if $3\le m\le 4$ and $2\le r\le 3$, or if $(m-1)(r-1)\le {r+1\choose2}$, then $\chi_{lat}(A(mK_n,K_r)-e) = n$.
\end{theorem}

\begin{proof} Without loss of generality, we may assume $e=v_{i,1}v_{i,r+1}$. Consider the local antimagic total labeling $f$ in the proof of Theorem~\ref{thm-A(mKnKr)}. For Case (1) with $m=4$ and $r=2$, we have for $A(4K_n,K_2)$, vertex $v_{i,1}$ and $v_{i,2}$ has degree $4n-7$, while $v_{i,j}$ has degree $n-1$ for $3\le j\le n$. We now have 
\begin{enumerate}[(i)]
  \item $w(v_{i,1}) - (4n-7)\ne w(v_{i,2}) - (4n-7)$;
  \item $w(v_{i,j}) - (4n-7)\ne w(v_{i,k}) - (n-1)$ for $j=1,2$ and $3\le k\le n$;
  \item $w(v_{i,j}) - (n-1)\ne w(v_{i,k}) - (n-1)$ for $3\le j<k\le n$.
\end{enumerate} 

\nt For all other cases, one can check that $f$ satisfies the conditions of Lemma~\ref{lem-lat-G-e}(ii) too. The details are omitted. Thus, the theorem holds.
\end{proof}

\nt By an argument similar to that for Theorem~\ref{thm-A(mKnKr)vK1}, we also have the the following theorem. The details are also omitted.

\begin{theorem}\label{thm-A(mKnKr)vK1-e} For $m\ge2$ and $n\ge r+1\ge 2$,  if $3\le m\le 4$ and $2\le r\le 3$, or if $(m-1)(r-1)\le {r+1\choose2}$, then $\chi_{la}((A(mK_n,K_r)-e)\vee K_1) = n+1$. \end{theorem} 

\begin{theorem}\label{thm-LAT-UmKn} For $m,n\ge 2$, $\chi_{lat}(A(mK_n,K_1)) = n$.   \end{theorem} 

\begin{proof} Note that $|V(A(mK_n,K_1))| + |E(A(mK_n,K_1))| = (mn-m+1)+m{n\choose 2} = m(n+2)(n-1)/2 + 1$. Let $u=v_{1,n}=\ldots =v_{m,n}$. Define $f: V(A(mK_n,K_1)) \cup E(A(mK_n,K_1)) \to [1,m(n+2)(n-1)/2+1]$ as follows.
\begin{enumerate}[(i)]
  \item $f(u) = m(n+2)(n-1)/2+1$,
  \item $f(v_{i,j})=m(n+2)(n-1)/2-i(n-1)+n-j$ for $1\le i\le m, 1\le j\le n-1$,
  \item $f(v_{i,j}v_{i,k})=m[n(j-1)-{j\choose2}+k-j-1]+i$ for $1\le i\le m, 1\le j < k\le n$.
\end{enumerate}

\ms\nt Observe that for $1\le i\le m, 1\le j\le n-1$, 
\begin{eqnarray} \label{eqn-w(vij)}
w(v_{i,j})& = & f(v_{i,j})+\sum^{j-1}_{k=1}f(v_{i,j}v_{i,k}) + \sum^{n}_{k=j+1}f(v_{i,j}v_{i,k})\nonumber\\
&= &  m(n+2)(n-1)/2-i(n-1)+n-j + \sum^{j-1}_{k=1}\Bigg(m\Big[n(k-1)-{k\choose2}+j-k-1\Big]+i\Bigg) + \nonumber\\ 
& &  \sum^{n}_{k=j+1} \Bigg(m\Big[n(j-1)-{j\choose2}+k-j-1\Big]+i\Bigg)\nonumber\\
&= &  m(n+2)(n-1)/2+n-j+m\Bigg[n{j-1\choose2}-{j\choose3}-{j\choose2}+(j-1)^2\Bigg] + \nonumber\\
& &  m(n-j)\Bigg[n(j-1)-{j\choose2}-j-1\Bigg] + m(n+j+1)(n-j)/2.
\end{eqnarray}

\ms\nt Moreover,
\begin{eqnarray}\label{eqn-w(u)}
w(u) &= & m(n+2)(n-1)/2+1+\sum^{m}_{i=1}\sum^{n-1}_{j=1} f(v_{i,j}v_{i,n}) \nonumber\\
& = & m(n+)(n-1)/2+1+ \sum^{m}_{i=1}\sum^{n-1}_{j=1} \Bigg(m\Big[n(j-1)-{j\choose2}+n-j-1\Big]+i\Bigg) \nonumber\\
& = & m(n+2)(n-1)/2+1+ m^2\Bigg(n{n\choose2} - (n-1)^2 -{n-1\choose3}\Bigg)-(n-1){m\choose2}. 
\end{eqnarray}

\ms\nt Therefore, $w(v_{1,j})=\cdots=w(v_{m,j})$ for $1\le j\le n-1$. Suppose $p=|V(A(mK_n,K_1))|$ and $q=|E(A(mK_n,K_1))|$. Thus, the labeling matrix, denoted $A(i,n,1)$, is as follow where $p+q$ is assigned once in $A(1,n,1)$.

\[\fontsize{10}{11}\selectfont\begin{array}{c|*{4}{c|}}
A(i,n,1) & v_{i,1} & v_{i,2} & v_{i,3} & v_{i,4} \\\hline
v_{i,1} & p+q-i(n-1)+n-2 & i & m+i & 2m+i  \\\hline
v_{i,2} & i & p+q-i(n-1)+n-3 & m(n-1)+i & mn+i \\\hline
v_{i,3} & m+i & m(n-1)+i & p+q-i(n-1)+n-4 & m(2n-3)+i \\\hline
v_{i,4} & 2m+i & mn+i & m(2n-3)+i & p+q-i(n-1)+n-5  \\\hline
\vdots & \vdots & \vdots & \vdots & \vdots  \\\hline
v_{i,j} & m(j-2)+i & m(n+j-4)+i & m(n(j-1)-{c\choose2}+2)+i & m(3n+j-11)+i  \\\hline
\vdots & \vdots & \vdots & \vdots & \vdots \\\hline
v_{i,n-1} & m(n-3)+i & m(2n-5)+i & m(2n-8)+i & m(4n-12)+i  \\\hline 
v_{i,n} & m(n-2)+i & m(2n-4)+i & m(3n-7)+i & m(4n-11)+i  \\\hline 
\end{array}\]

\[\fontsize{10}{11}\selectfont\begin{array}{c|*{5}{c|}}
 &  \ldots & v_{i,j} & \ldots & v_{i,n-1} & v_{i,n}\\\hline
v_{i,1} & \ldots & m(j-2)+i & \ldots & m(n-3)+i & m(n-2)+i \\\hline
v_{i,2} & \ldots & m(n+j-4)+i & \ldots & m(2n-5)+i & m(2n-4)+i\\\hline
v_{i,3} & \ldots & m(2n+j-7)+i & \ldots & m(2n-8)+i & m(3n-7)+i\\\hline
v_{i,4} & \ldots & m(3n+j-11)+i & \ldots & m(4n-12)+i & m(4n-11)+i \\\hline
\vdots & \ddots & \vdots & \vdots & \vdots & \vdots \\\hline
v_{i,j} & \ldots & p+q-i(n-1)+n-j+1 & \ldots & m(nj-{j\choose2}-j-2)+i & m(nj-{j\choose2}-j)+i \\\hline
\vdots & \vdots & \vdots & \ddots & \vdots & \vdots \\\hline
v_{i,n-1} & \ldots & m(nj-{j\choose2}-j-2)+i & \ldots & p+q-i(n-1) & m(n-2)(n+1)/2+i  \\\hline 
v_{i,n} & \ldots & m(nj-{j\choose2}-j)+i & \ldots & m(n-2)(n+1)+i & \star  \\\hline 
\end{array}\]

\ms\nt Obviously, for $1\le j \le n-1$, (i) there is exactly one identical entry in rows $j$ and $j+1$ given by $f(v_{i,j+1}v_{i,j})=f(v_{i,j}v_{i,j+1})$, (ii) $f(v_{j+1,j+1})+1=f(v_{j,j})$ for $j+1<n$, (iii) all other entries in column $k$ of row $j+1$ are larger than the corresponding entries in row $j$ for $k\ne j,j+1$, and (iv) $f(u)$ in row $n$ is the largest of all labels. Thus, $w(v_{i,1}) < w(v_{i,2}) < \cdots < w(v_{i,n})$.   Therefore, $f$ is a local antimagic total labeling that induces $n$ distinct weights. Since $\chi(A(mK_n,K_1))=n$, the theorem holds. 
\end{proof} 

\begin{example} For $m=2,n=5$, the matrices $A(i,5,1), i=1,2,$ are given below.\\

\nt\begin{minipage}[b]{8.5cm}
$\fontsize{9}{10}\selectfont\begin{array}{c||*{5}{c|}|c}
A(1,5,1) & v_{1,1} & v_{1,2} & v_{1,3} & v_{1,4} & v_{1,5} & w(v_{i,j})\\\hline\hline
v_{1,1} & 28 & 1 & 3 & 5 & 7 & 44\\\hline
v_{1,2} & 1 & 27 & 9 & 11 & 13 & 61\\\hline
v_{1,3} & 3 & 9 & 26 & 15 & 17 & 70\\\hline
v_{1,4} & 5 & 11 & 15 & 25 & 19 & 75\\\hline
u=v_{1,5} & 7 & 13 & 17 & 19 & 29 & \star\\\hline
\end{array}$\\
$f(u) + \sum^4_{k=1} f(v_{1,5}v_{1,k})=85$\\ $w(u)=145$.
\end{minipage}
\begin{minipage}[b]{8.5cm}
$\fontsize{9}{10}\selectfont\begin{array}{c||*{5}{c|}|c}
A(2,5,1) & v_{2,1} & v_{2,2} & v_{2,3} & v_{2,4} & v_{2,5} & w(v_{i,j})\\\hline\hline
v_{2,1} & 24 & 2 & 4 & 6 & 8 & 44\\\hline
v_{2,2} & 2 & 23 & 10 & 12 & 14 & 61\\\hline
v_{2,3} & 4 & 10 & 22 & 16 & 18 & 70\\\hline
v_{2,4} & 6 & 12 & 16 & 21 & 20 & 75\\\hline
u=v_{1,5} & 8 & 14 & 18 & 20 & \star & \star\\\hline
\end{array}$\\
$\sum^4_{k=1} f(v_{1,5}v_{1,k})=60$\\
\end{minipage}\\
\end{example}

\begin{theorem}\label{thm-LAT-UmKn-e} Let $e$ be an edge of 
\begin{enumerate}[(I)]
  \item $A(mK_2,K_1)$, $m\ge 2$, then $\chi_{lat}(A(mK_2,K_1)-e)=2$;
  \item $A(mK_n,K_1)$, $m\ge 2, n\ge 3$, that is not incident to the vertex $u=v_{1,n}=\cdots=v_{m,n}$, then $\chi_{lat}(A(mK_n,K_1)-e) = n$,
  \item $A(mK_n,K_1)$, $m,n\ge 3$, that is incident to the vertex $u=v_{1,1}=\cdots=v_{m,1}$, then $\chi_{lat}(A(mK_n,K_1)-e) = n$. 
 \end{enumerate} 
\end{theorem} 

\begin{proof} (I). Note that $A(mK_2,K_1)-e\cong K_{1,m-1} + K_1$. Let $V(K_{1,m-1}+K_1) = \{x,u,v_i\,|\,1\le i\le m-1\}$ and $E(K_{1,m-1}+K_1) = \{uv_i\,|\,1\le i\le m-1\}$. For $1\le i\le m-1$, define $f: V(K_{1,m-1}+K_1)\cup E(K_{1,m-1}+K_1)\to [1,2m]$ such that $f(x)=2m$, $f(v_i)=i$, $f(uv_i)=2m-i$, $f(u)= m$. Thus, $w(x)=w(v_i)=2m\ne w(u)=m(3m-1)/2$. Thus, $\chi_{lat}(A(mK_2,K_1)-e) = 2$.

\ms\nt (II). Suppose $e=v_{1,1}v_{1,2}$ not incident to $u=v_{1,n}=\cdots=v_{m,n}$. Let $f$ be the local antimagic total labeling of $A(mK_n,K_1)$ in Theorem~\ref{thm-LAT-UmKn}. Define $h: V(A(mK_n,K_1)-e)\cup E(A(mK_n,K_1)-e)\to [1,m(n+2)(n-1)/2]$ such that $g = f - 1$ for every vertex and edge of $A(mK_n,K_1)-e$. Thus, $w_g(v_{i,j}) = w_f(v_{i,j}) - n$ for $1\le j\le n-1$, and $w_g(u)=w_f(u)-m(n-1)-1$. Note that $f(v_{1,1}v_{1,2})=0$ means edge $e=v_{1,1}v_{1,2}$ is deleted. Clearly, $w_g(x) = w_g(y)$ if and only if $w_f(x) = w_f(y)$ for $x,y\ne u$. Since $w_g(v_{i,n-1}) > w_g(v_{i,j})$ for $j<n-1$, by Lemma~\ref{lem-lat-G-e}, suffice to show that $A=w_f(u) - w_f(v_{1,n-1}) - (m-1)(n-1) > 0$. From Equations~\ref{eqn-w(u)} and~\ref{eqn-w(vij)}, we have

\begin{eqnarray*} 
A& = & \Bigg\{m(n+2)(n-1)/2+1+ m^2\Bigg(n{n\choose2} - (n-1)^2 -{n-1\choose3}\Bigg)-(n-1){m\choose2}\Bigg\}- \\
& & \Bigg\{m(n+2)(n-1)/2+n-j+m\Bigg[n{j-1\choose2}-{j\choose3}-{j\choose2}+(j-1)^2\Bigg] + \\
& &  m(n-j)\Bigg[n(j-1)-{j\choose2}-j-1\Bigg] + m(n+j+1)(n-j)/2\Bigg\} - (m-1)(n-1)\\
& = & \frac{m^2n^3}{3} - \frac{mn^3}{3} - \frac{m^2n^2}{2}+\frac{mn^2}{2}-\frac{m^2n}{3}+\frac{4mn}{3}+n+\frac{m^2}{2}-\frac{5m}{2}-1.\\
\end{eqnarray*}

\nt Clearly, $A>0$ if $m=2$ or $n=3$. Suppose $m\ge 3$ and $n\ge 4$, then 
\begin{eqnarray*}
A&=& \frac{1}{6}(2m^2n^3-2mn^3-3m^2n^2+3mn^2-2m^2n) + \frac{4mn}{3}+n+\frac{m^2}{2}-\frac{5m}{2}-1\\
& > &\frac{1}{6}(m^2n^3-2mn^3+m^2n^3-3m^2n^2+3mn^2-2m^2n) \\
& \ge &\frac{1}{6}(4m^2n^2-3m^2n^2+3mn^2-2m^2n) > 0.
\end{eqnarray*}

\nt Thus, $g$ is a local antimagic total labeling that induces $n$ distinct weights. Since $\chi(A(mK_n,K_1)-e)=n$, we have $\chi_{lat}(A(mK_n,K_1)-e) = n$.

\ms\nt (III). Suppose $e=v_{1,1}v_{1,2}$ incident to $u=v_{1,1}=\cdots = v_{m,1}$ so that $e$ is incident to $u$. Define $f: V(A(mK_n,K_1)) \cup E(A(mK_n,K_1)) \to [1,m(n+2)(n-1)/2]$ as follows.

\begin{enumerate}[(i)]
  \item $f(u) = m(n+2)(n-1)/2$,
  \item $f(v_{i,j})=m(n+2)(n-1)/2-i(n-1)+n-j$ for $1\le i\le m, 2\le j\le n$,
  \item $f(v_{i,j}v_{i,k})=m[{k-1\choose2}+j-1]+i-1$ for $1\le i\le m, 1\le j < k\le n$.
\end{enumerate}

\nt Note that $f(v_{1,1}v_{1,2})=0$ means edge $e=v_{1,1}v_{1,2}$ is deleted. Observe that for $1\le i\le m, 2\le j\le n$, 

\begin{eqnarray} \label{eqn-w(vij)-2}
w(v_{i,j})& = & f(v_{i,j})+\sum^{j-1}_{k=1}f(v_{i,j}v_{i,k}) + \sum^{n}_{k=j+1}f(v_{i,j}v_{i,k})\nonumber\\
&= &  m(n+2)(n-1)/2-i(n-1)+n-j + \sum^{j-1}_{k=1}\Bigg(m\Bigg[{j-1\choose2}+k-1\Bigg]+i-1\Bigg) + \nonumber\\ 
& &  \sum^{n}_{k=j+1} \Bigg(m\Bigg[{k-1\choose2}+j-1\Bigg]+i-1\Bigg)\nonumber\\
&= &  m(n+2)(n-1)/2-j+1 + m(j-1)\Big[{j-1\choose2}+n-j-1\Big]+\nonumber\\
& & m{j\choose2}+m{n\choose3}-m{j\choose3}.
\end{eqnarray}

\nt Moreover,
\begin{eqnarray}\label{eqn-w(u)-2}
w(u) &= & m(n+2)(n-1)/2+\sum^{m}_{i=1}\sum^{n}_{k=2} f(v_{i,1}v_{i,k}) \nonumber\\
& = & m(n+)(n-1)/2+ \sum^{m}_{i=1}\sum^{n}_{k=2} \Bigg(m{k-1\choose2}+i-1\Bigg) \nonumber\\
& = & m(n+2)(n-1)/2+ m^2{n\choose3} + (n-1){m\choose2}. 
\end{eqnarray}

\nt Listing the labels in matrix form as in the proof of Theorem~\ref{thm-LAT-UmKn}, we have $w(v_{i,2}) < \cdots < w(v_{i,n})$. Suffice to show that $w(u) > w(v_{1,n})$. Since $m,n\ge 3$, we have
\begin{eqnarray*}
w(u)-w(v_{1,n})&=& m^2{n\choose3} + (n-1){m\choose2} + n-1- m(n-1)\Big[{n-1\choose2}-1\Big]-m{n\choose2}\\
&=&\frac{m^2n^3}{6}-\frac{mn^3}{2}-\frac{m^2n^2}{2}+\frac{3mn^2}{2}+\frac{5m^2n}{6}-\frac{3mn}{2}+n-\frac{m^2}{2}+\frac{m}{2}-1\\
&>&0
\end{eqnarray*}
Thus, $f$ is a local antimagic total labeling that induces $n+1$ distinct weights. Since $\chi(A(mK_n,K_1)-e)=n$, we have $\chi_{lat}(A(mK_n,K_1)-e) = n$. 
\end{proof}

\begin{example} For Theorem~\ref{thm-LAT-UmKn-e}(III) with $m=3, n=4$ and $e=v_{1,1}v_{1,2}$ incident to $u=v_{1,1}=\cdots = v_{2,1}$, the matrices $A(i,4,1)$, $i=1,2,3$, are given below. Note that the label of the common vertex only appears in $A(1,4,1)$.\\

\hskip1.3cm\begin{minipage}[b]{7cm}
$\fontsize{9}{10}\selectfont\begin{array}{c||*{4}{c|}|c}
A(1,4,1) & v_{1,1} & v_{1,2} & v_{1,3} & v_{1,4} & w(v_{1,j})\\\hline\hline
u=v_{1,1} & 27 & 0 & 3 & 9 & \star\\\hline
v_{1,2} & 0 & 26 & 6 & 12 & 53\\\hline
v_{1,3} & 2 & 4 & 25 & 15 & 58\\\hline
v_{1,4} & 6 & 8 & 10 & 24 & 67\\\hline
\end{array}$\\
\hskip2cm $f(u) + \sum^4_{k=3} f(v_{1,1}v_{1,k})=39$\\\hskip2cm $w(u)=72$.
\end{minipage}
\begin{minipage}[b]{8cm}
$\fontsize{9}{10}\selectfont\begin{array}{c||*{4}{c|}|c}
A(2,4,1) & v_{2,1} & v_{2,2} & v_{2,3} & v_{2,4} & w(v_{2,j})\\\hline\hline
u=v_{2,1} & \star & 1 & 4 & 10 & \star\\\hline
v_{2,2} & 1 & 23 & 7 & 13 & 53\\\hline
v_{2,3} & 3 & 5 & 22 & 16 & 58\\\hline
v_{2,4} & 7 & 9 & 11 & 21 & 67\\\hline
\end{array}$\\
$\sum^4_{k=2} f(v_{2,1}v_{2,k})=15$\\
\end{minipage}\\
\begin{table}[h]\hskip2cm
$\fontsize{9}{10}\selectfont\begin{array}{c||*{4}{c|}|c}
A(3,4,1) & v_{2,1} & v_{2,2} & v_{2,3} & v_{2,4} & w(v_{3,j})\\\hline\hline
u=v_{3,1} &\star & 2 & 5 & 11 & \star\\\hline
v_{3,2} & 1 & 20 & 8 & 14 & 53\\\hline
v_{3,3} & 3 & 5 & 19 & 17 & 58\\\hline
v_{3,4} & 7 & 9 & 11 & 18 & 67\\\hline
\end{array}$
\hskip3mm $\sum^4_{k=2} f(v_{3,1}v_{3,k})=18$
\end{table}
\end{example}

\nt Suppose $m=2$. By computer search, we found that the labeling defined for Theorem~\ref{thm-LAT-UmKn-e}(III) gives $w(u)\ne w(v_{1,j})$ for $2\le j\le n\le 10000$.

\begin{theorem}\label{thm-LAT-UmKn-v-K1} Let $e$ be an edge of 
\begin{enumerate}[(a)]
  \item $A(mK_n,K_1)$, $m,n\ge 2$, then $\chi_{la}(A(mK_n,K_1)\vee K_1)=n+1$.
  \item $A(mK_2,K_1)$, $m\ge 2$, then $\chi_{la}((A(mK_2,K_1)-e)\vee K_1)=3$;
  \item $A(mK_n,K_1)$, $m\ge 2, n\ge 3$, that is not incident to the vertex $u=v_{1,n}=\cdots=v_{m,n}$, then $\chi_{la}((A(mK_n,K_1)-e)\vee K_1) = n+1$,
  \item $A(mK_n,K_1)$, $m,n\ge 3$, that is incident to the vertex $u=v_{1,1}=\cdots=v_{m,1}$, then $\chi_{la}((A(mK_n,K_1)-e)\vee K_1) = n+1$. 
 \end{enumerate} 
 \end{theorem}

\begin{proof} (a) Let $G=A(mK_n,K_1)$, $m\ge 2, n\ge 3$. Consider $G$ and the local antimagic total labeling $f$ of $G$ as in the proof of Theorem~\ref{thm-LAT-UmKn}. We note that $f(v_{1,n})$ is the largest assigned label. Moreover, in each $M_i$, $1\le i\le m$ and $1\le k\le n-1$, $f(u_{i,k}) > f(v_{i,n}v_{i,k})$. Thus, sum of all the vertex labels under $f$ is larger than all the induced vertex weights under $f$. Since $\chi_{lat}(G)) = \chi(G\vee K_1)-1=(n+1)-1$, by Theorem~\ref{thm-K1VG}(b), we conclude that $\chi_{la}(A(mK_n,K_1) \vee K_1)=n+1$.

\ms\nt (b) Let $G=A(mK_2,K_1)-e$, $m\ge 2$. Consider $G$ and the local antimagic total labeling $f$ of $G$ as in the proof of Theorem~\ref{thm-LAT-UmKn-e}(I). It is easy to show that $\chi_{la}(A(mK_2,K_1)-e)=3$ for $m=2,3$. If $m\ge 4$, we note that the sum of all vertex labels under $f$ is $m(m+5)/2$ which is not equal to all the induced vertex weights under $f$. Since $\chi_{lat}(G) = \chi(G\vee K_1)-1=3-1$, by Theorem~\ref{thm-K1VG}(b), we conclude that $\chi_{la}((A(mK_n,K_1)-e) \vee K_1)=3$.

\ms\nt (c) Let $G=A(mK_n,K_1)-e$, $m\ge 2, n\ge 3$, where $e=v_{1,1}v_{1,2}$ is not incident to the vertex $u=v_{1,n}=\cdots=v_{m,n}$. Consider $G$ and the local antimagic total labeling $f$ of $G$ as in the proof of Theorem~\ref{thm-LAT-UmKn-e}(II). It is easy to show that the sum of all vertex labels under $f$ is $\frac{1}{2}m(n-1)(n+1)[m(n-1)+1]$ which is larger than all the induced vertex weights under $f$. Since $\chi_{lat}(G) = \chi(G\vee K_1)-1=n-1$, by Theorem~\ref{thm-K1VG}(b), we conclude that $\chi_{la}((A(mK_n,K_1)-e) \vee K_1)=n+1$.   

\ms\nt (d) Let $G=A(mK_n,K_1)-e$, $m,n\ge 3$, where $e=v_{1,1}v_{1,2}$ is incident to the vertex $u=v_{1,1}=\cdots=v_{m,1}$. Consider $G$ and the local antimagic total labeling $f$ of $G$ as in the proof of Theorem~\ref{thm-LAT-UmKn-e}(III). It is easy to show that the sum of all vertex labels under $f$ is $\frac{1}{2}m(n-1)(n+1)[m(n-1)+1]$ which is larger than all the induced vertex weights under $f$. Since $\chi_{lat}(G) = \chi(G\vee K_1)-1=n-1$, by Theorem~\ref{thm-K1VG}(b), we conclude that $\chi_{la}((A(mK_n,K_1)-e) \vee K_1)=n+1$.  
\end{proof}


\begin{theorem}\label{thm-lat-1reg} For $m\ge 2$ and even $n\ge 2$, $\chi_{lat}(mK_n)=\chi_{lat}(mK_n - e)= n$. Otherwise, $n\le \chi_{lat}(mK_n)\le m+n-1$. \end{theorem}

\begin{proof} We first note that for $G=mK_n$, $\chi_{lat}(G), \chi_{lat}(G-e)\ge n$, and that $\chi_{la}(G \vee K_1), \chi_{la}((G\vee K_1) - e) \ge n+1$. For $1\le i\le m$, let $V(mK_2) = \{u_i,v_i\}$ and $E(mK_2) = \{u_iv_i\}$. Define $f: V(mK_2)\cup E(mK_2) \to [1,3m]$ such that $f(u_iv_i)=i$, $f(u_i) = 2n+1-i$ and $f(v_i) = 3n+1-i$. Thus, $f$ is a local antimagic total labeling that induces 2 distinct vertex weights.

\ms\nt Consider even $n\ge4$. Clearly, $|V(mK_n)|=mn$ and $|E(mK_n)| = mn(n+1)/2$. Let $[t] = \{3(t-1)+1, 3(t-1)+2, \ldots, 3(t-1)+m\}$. Let $\mathcal M$ be a symmetric $n\times n$ matrix. Let $M_i$ be the labeling matrix of the $i$-th $K_n$. We define a total labeling of the $i$-th $K_n$ by assigning appropriate integers in $[1,mn(n+1)/2]$ as follows. 

\begin{enumerate}[(i)]
  \item For matrix $\mathcal M$, assign $[1]$ to $[n(n+1)/2]$ beginning with the main diagonal entry horizontally to the right and consecutively from row $1$ to row $n$ such that column $k$ is the transpose of row $k$.
  \item Beginning with the (1,1)-entry, bold the entries alternately for every row and column.
  \item For the $i$-th labeling matrix, $M_i$, the entry that corresponds to a bold $[t]$ is the $i$-th integer in $[t]$ whereas the entry that corresponds to a non-bold $[t]$ is the $(m+1-i)$-th integer in $[t]$.
  \item For each row, observe that the sum of the first, second, $\dots, m/2$-th pair of entries are equal respectively. Thus, $M_1$ to $M_m$ always have equal row sum for row 1 to row $n$ respectively.
  \item Since the numbers are assigned consecutively from row 1 to row $n$, the row sums form an increasing sequence. 
\end{enumerate} 

\nt Consequently, the labeling matrices represent a local antimagic total labeling of the $mK_n$ that induces $n$ distinct vertex weights. Thus, $\chi_{lat}(mK_n)\le n$. 

\ms\nt Observe that the largest row sum is $R=mn+\sum^{n}_{j=2}(m(jn-{j\choose2})-m+1) = m[n{n+1\choose2}-{n+1\choose3}] - n(m-1)/2$, and the sum of all vertex labels is $S=\sum^n_{j=1}(m^2[(j-1)n-{j-1\choose2}]+{m+1\choose2}) = 2m^2{n+1\choose3}+n{m+1\choose2}$. Thus, $S-R = m^2n^3/3-mn^3/2-mn^2/2+m^2n/6+mn-n/2+m^4/6-m^2/6 > 0$. By Theorem~\ref{thm-K1VG}(b), we have $\chi_{la}(mK_n \vee K_1)\le n+1$.  

\ms\nt We now consider $mK_n - e$. If $n=2$, the theorem clearly holds by Lemma~\ref{lem-regular}. We now modify the labeling of $mK_n$ by swapping the vertex label $i$ and the incident edge label $2m+i$ for each $1\le i\le m$. Thus, all vertex labels remained unchanged except that the third vertex label of each $K_n$ is now reduced by $2m$. Since the original second row sum is less than the original third row sum by more than $2m$, we now have a local antimagic total labeling of $mK_n$ that induces $n$ distinct weight. Since the edge joining the first and third vertex of the first $K_n$ is labeled by $1$, by Lemmas~\ref{lem-non-regular} and~\ref{lem-lat-G-e}, we conclude that $\chi_{lat}(mK_n - e)\le n$.    



\ms\nt Consider odd $n\ge 3$. By modifying the labeling function for $A(mK_n,K_1)$ in the proof of Theorem~\ref{thm-LAT-UmKn}, it is easy to get a local antimagic total labeling for $mK_n$, $m\ge 2, n\ge 3$, that induces $m+n-1$ distinct weights. The theorem holds. 
\end{proof}

\nt We note that $mK_n \vee K_1\cong A(mK_{n+1},K_1)$ as in Theorem\ref{thm-A(mKnKr)vK1}.

\begin{example} We take $m=3$ and $n=6$. The labeling matrices of all the cases with exact $\chi_{lat}$ and $\chi_{la}$ are given below.

\nt\begin{minipage}[b]{9cm}
$\fontsize{9}{10}\selectfont\begin{array}{c||*{6}{c|}}
\mathcal M & v_{1,1} & v_{1,2} & v_{1,3} & v_{1,4} & v_{1,5} & v_{1,6}  \\\hline\hline
v_{1,1} & {\bf[1]} & [2] & {\bf[3]} & [4] & {\bf[5]} & [6]  \\\hline
v_{1,2} & [2] & {\bf[7]} & [8] & {\bf[9]} & [10] & {\bf[11]} \\\hline
v_{1,3} & {\bf[3]} & [8] & {\bf[12]} & [13] & {\bf[14]} & [15]   \\\hline
v_{1,4} & [4] & {\bf[9]} & [13] & {\bf[16]} & [17] & {\bf[18]}   \\\hline
v_{1,5} & {\bf[5]} & [10] & {\bf[14]} & [17] & {\bf[19]} & [20]  \\\hline
v_{1,6} & [6] & {\bf[11]} & [15] & {\bf[18]} & [20] & {\bf[21]}   \\\hline
\end{array}$\\
\end{minipage}
\begin{minipage}[b]{9cm}
$\fontsize{9}{10}\selectfont\begin{array}{c||*{6}{c|}|c}
M_1 & v_{2,1} & v_{2,2} & v_{2,3} & v_{2,4} & v_{2,5} & v_{2,6} &  \mbox{sum}\\\hline\hline
v_{2,1} & 1 & 6 & 7 & 12 & 13 & 18 & 57  \\\hline
v_{2,2} & 6 & 19 & 24 & 25 & 30 & 31 & 135 \\\hline
v_{2,3} & 7 & 24 & 36 & 37 & 42 & 43 & 189  \\\hline
v_{2,4} & 12 & 25 & 37 & 48 & 49 & 54 & 225  \\\hline
v_{2,5} & 13 & 30 & 42 & 49 & 55 & 60 & 249  \\\hline
v_{2,6} & 18 & 31 & 43 & 54 & 60 & 61 & 267  \\\hline
\end{array}$\\
\end{minipage}\\\\
\nt\begin{minipage}[b]{9cm}
$\fontsize{9}{10}\selectfont\begin{array}{c||*{6}{c|}|c}
M_2 & v_{3,1} & v_{3,2} & v_{3,3} & v_{3,4} & v_{3,5} & v_{3,6} & \mbox{sum}\\\hline\hline
v_{3,1} & 2 & 5 & 8 & 11 & 14 & 17 & 57  \\\hline
v_{3,2} & 5 & 20 & 23 & 26 & 29 & 32 & 135  \\\hline
v_{3,3} & 8 & 23 & 35 & 38 & 41 & 44 & 189  \\\hline
v_{3,4} & 11 & 26 & 38 & 47 & 50 & 53 & 225  \\\hline
v_{3,5} & 14 & 29 & 41 & 50 & 56 & 59 & 249  \\\hline
v_{3,6} & 17 & 32 & 44 & 53 & 59 & 62 & 267  \\\hline
\end{array}$\\
\end{minipage}
\begin{minipage}[b]{9cm}
$\fontsize{9}{10}\selectfont\begin{array}{c||*{6}{c|}|c}
M_3 &v_{4,1} & v_{4,2} & v_{4,3} & v_{4,4} & v_{4,5} & v_{4,6} & \mbox{sum}\\\hline\hline
v_{4,1} & 3 & 4 & 9 & 10 & 15 & 16 & 57  \\\hline
v_{4,2} & 4 & 21 & 22 & 27 & 28 & 33 & 135  \\\hline
v_{4,3} & 9 & 22 & 36 & 37 & 42 & 43 & 189  \\\hline
v_{4,4} & 10 & 27 & 37 & 48 & 49 & 54 & 225  \\\hline
v_{4,5} & 15 & 28 & 42 & 49 & 57 & 58 & 249  \\\hline
v_{4,6} & 16 & 33 & 43 & 54 & 48 & 63 & 267  \\\hline
\end{array}$\\
\end{minipage}\\
\nt The above matrices give $\chi_{lat}(3K_6) = 6$. Let $v$ be the vertex of $K_1$. If the main diagonal labels are the edge labels of $3K_6 \vee K_1$, then $v$ has label $666$. Thus, the matrices give $\chi_{la}(3K_6 \vee K_1)=7$. Deleting edge $vv_{1,1}$ of label 1 of $3K_3\vee K_1$ and reducing all edge labels by 1, we get a local antimagic labeling of $(3K_6\vee K_1)-vv_{1,1}$. Thus, by symmetry, $\chi_{la}((3K_6\vee K_1)-e)=7$ for $e$ not belong to any $K_6$. If we swap the labels of 1 and 7, 2 and 8, and 3 and 9, delete the edge $v_{1,1}v_{1,3}$ that has label 1 and reduce all other labels by 1, we get $\chi_{lat}(3K_6 - e)=6$. Now, if the main diagonal labels are the edge labels of $3K_6 \vee K_1$, then we have $\chi_{la}(3K_6\vee K_1) - e)=7$ for $e$ that belongs to any $K_6$.
\end{example}

\nt Let $v$ be the central vertex of the wheel $W_n = C_n \vee K_1$ of order $n+1\ge 4$. Denote by $A(mW_n,v)$ the amalgamation of $m\ge 1$ copies of $W_n$ at vertex $v$ where $A(W_n,v) = W_n$. In~\cite{Slamin+NMDK}, the authors claimed that $\chi_{lsat}(A(mW_n,v)) = 3 $ (or 4) for even (or odd) $n$. However, the labeling function defined in the proof contains errors that can be corrected as follow. 

\begin{theorem}\label{thm-A(mWn,v)} For $m\ge 1, n\ge 3$, $$\chi_{lat}(A(mW_n,v)) = \begin{cases} 3  & \mbox{ if } m \mbox{ is even,} \\ 4 &\mbox{ if } m \mbox{ is odd.} \end{cases}$$
\end{theorem}

\begin{proof} Let $G = A(mW_n,v)$, $V(G) = \{v, v_{i,j}\,|\,1\le i\le m, 1\le j\le n\}$ and $E(G) = \{vv_{i,j}\,|\,1\le i\le m, 1\le j\le n\} \cup \{v_{i,1}v_{i,n}, v_{i,j}v_{i,j+1}\,|\,1\le i\le m, 1\le j\le n-1\}$. Suppose $n$ is odd. Define a total labeling $f: V(G)\cup E(G) \to [1,3mn+1]$ as follows.
\begin{enumerate}[(i)]
  \item $f(v)=1$ and $f(v_{i,1}) = m(n-1) + i+1$,
  \item $f(v_{i,j}) = m(j-2)+i+1$ for $1\le i\le m, 2\le j\le n$,
  \item $f(vv_{i,j}) = mn+1+m(j-3)+i$ for $1\le i\le m$ and odd $j=3,5,\ldots,n-1$,
  \item $f(vv_{i,j}) = m(n-1)+1+mj+i$ for $1\le i\le m$ and even $j=2,4,\ldots,n-2$,
  \item $f(vv_{i,1}) = m(2n-1)+1+i$ for $1\le i\le m$,
  \item $f(v_{i,1}v_{i,n}) = 2mn + m + 2 - i$ for $1\le i\le m$,
  \item $f(v_{i,j}v_{i,j+1}) = 3mn + m(1-j)-i+2$ for $1\le i\le m, 1\le j\le n-1$. 
\end{enumerate}
\nt It is easy to check that $w(v)=1+3{mn+1\choose2}$, $w(v_{i,1}) = 8mn-m+6$, $w(v_{i,j}) = 7mn-2m+6$ for odd $j > 1$, and $w(v_{i,j}) = 7mn+6$ for even $j$. Thus, $f$ is a local (super) antimagic total labeling of $G$. Since $4=\chi(G)\le \chi_{lat}(G)\le 4$, the theorem holds.

\ms\nt Suppose $n$ is even. Define a total labeling $f: V(G)\cup E(G) \to [1,3mn+1]$ such that the labels are as for odd $n$ for all the vertices and the edges $v_{i,1}v_{i,n}$. Moreover,
\begin{enumerate}[(i)]
  \item $f(vv_{i,j}) = 2mn-mj+i+1$ for $1\le i\le m$, $1\le j\le n$,
  \item $f(v_{i,j}v_{i,j+1}) = 2mn+m(j+1)+2-i$ for $1\le i\le m$, $1\le j\le n-1$.
\end{enumerate} 
\nt It is easy to check that $w(v) = 1+3{mn+1\choose2}$, $w(v_{i,j}) = 7mn + m + 6$ for odd $j$, and $w(v_{i,j}) = 7mn - m+6$ for even $j$. Thus, $f$ is a local (super) antimagic total labeling of $G$. Since $3=\chi(G)\le \chi_{lat}(G)\le 3$, the theorem holds.
\end{proof} 





\nt Note that $W_n - e$ is the fan graph $F_n$ if $e$ is an edge not incident to the vertex $v$. 

\begin{corollary}\label{cor-Fn} For $n\ge 3$, $\chi_{lat}(F_n) =  3$ if $n$ is even, and $3\le \chi_{lat}(F_n)\le 4$, otherwise.  \end{corollary} 

\nt Note that the above theorem also holds for $\chi_{lsat}(A(mK_n,v))$. By Theorem~\ref{thm-K1VG}, Lemma~\ref{lem-non-regular} or Lemma~\ref{lem-lat-G-e}, we also have

\begin{corollary}\label{} For $n\ge 3$ and $e$ an edge not incident to vertex $v$, $$\chi_{lat}(A(mW_n,v)-e) = \begin{cases} 3 \mbox{ if } n \mbox{ is even,} \\ 4 \mbox{ if } n \mbox{ is odd.}\end{cases}$$ \end{corollary}

\begin{corollary}\label{cor-A(mWn,v)VK1} For $n\ge 3$ and $e$ an edge of $A(mW_n,v)$, $$\chi_{la}(A(mW_n,v) \vee K_1) = \chi_{la}((A(mW_n,v) \vee K_1) - e) = \begin{cases} 4 \mbox{ if } n \mbox{ is even,} \\ 5 \mbox{ if } n \mbox{ is odd.}\end{cases}$$ \end{corollary}

\nt We now determine $\chi_{lat}(K_m \odot K_n)$ for $m\ge 1,n\ge 2$.

\begin{theorem}\label{thm-Km-odot-Kn} For $n\ge 2$, $\chi_{lat}(K_1\odot K_n) = n+1$ and $\max\{m,n+1\}\le \chi_{lat}(K_m\odot K_n)\le m+n$ if $m\ge 2$. \end{theorem}

\begin{proof} Let $G = K_m \odot K_n$, $m,n\ge 2$. Note that we can view $G$ as a graph obtained from $m$ copies of $K_{n+1}$ by taking a vertex of each $K_{n+1}$, and join these vertices by edges pairwise to form a $K_m$. Thus, we may have $V(G) = \{v_{i,j}\,|\,1\le i\le m, 1\le j\le n+1\}$ and $E(G) = \{v_{i,j}v_{i,k}\,|\,1\le j < k \le n+1\} \cup \{v_{i,n+1}v_{j,n+1}\,|\,1\le i<j\le m\}$. Note that vertices $v_{i,n+1}$, $1\le i\le m$, form the $K_m$ subgraph with $|V(G)|=m(n+1)$ and $|E(G)| = m{n+1\choose2}+{m\choose2}$.

\ms\nt Define $f: V(G)\cup E(G)\to [1,m(n+1)+m{n+1\choose2}+{m\choose2}]$ as follows.
\begin{enumerate}[(i)]
  \item $f(v_{i,j}v_{i,k})=m[{k-1\choose2}+j-1]+i$ for $1\le i\le m, 1\le j < k\le n+1$;
  \item $f(v_{i,j}) = m{n+1\choose2}+(m-i+1)n - j+1$ for $1\le j\le n$;
  \item $f(v_{i,n+1}) = m{n+1\choose2} + mn+i$ for $1\le i\le m$;
  \item $f(v_{i,n+1}v_{k,n+1}) = m{n+1\choose2}+m(n+1) + {k-1\choose2}+i$ for $1\le i < k\le m$.
\end{enumerate}

\nt Observe that for $1\le i\le m$, 
\begin{eqnarray}\label{eqn-w(vi1)-odot}
w(v_{i,n+1})&=&f(v_{i,n+1}) + \sum^{i-1}_{k=1} f(v_{i,n+1}v_{k,n+1}) + \sum^m_{k=i+1} f(v_{i,n+1}v_{k,n+1}) + \sum^{n}_{j=1} f(v_{i,j}v_{i,n+1}) \nonumber\\
&=& m{n+1\choose2} + mn+i + \sum^{i-1}_{k=1} \Bigg(m{n+1\choose2}+m(n+1) + {i-1\choose2}+k\Bigg) + \nonumber\\
& & \sum^m_{k=i+1} \Bigg(m{n+1\choose2}+m(n+1) + {k-1\choose2}+i\Bigg) + \sum^{n}_{j=1}\Bigg[m{n\choose2}+m(j-1)+i\Bigg] \nonumber\\
&=& m{n+1\choose2} + mn+i + (i-1)\Bigg(m{n+1\choose2}+m(n+1) + {i-1\choose2}\Bigg)+{i\choose2} + \nonumber\\
& & (m-i)\Bigg(m{n+1\choose2}+m(n+1)+i\Bigg)  + {m\choose3}+mn{n\choose2}+m{n\choose2}+ni\nonumber\\
&=& m{n+1\choose2} + n(m+i)+i + (m-1)\Bigg(m{n+1\choose2}+m(n+1)\Bigg)+(i-1){i-1\choose2}+\nonumber\\
& & {i\choose2} + (m-i)i+ {m\choose3}+m(n+1){n\choose2}\nonumber\\
&=& m{n+1\choose2} + n(m+i)+  (m-1)\Bigg(m{n+1\choose2}+m(n+1)\Bigg)+mi+ {m\choose3}+m(n+1){n\choose2}+\nonumber\\
& & i^3/2-5i^2/2+3i-1
\end{eqnarray} 
which is an increasing function. 

\nt Moreover, for $1\le i\le m$, $1\le j \le n$, 
\begin{eqnarray}\label{eqn-w(vij-odot)}
w(v_{i,j})&=& f(v_{i,j}) + \sum^{j-1}_{k=1} f(v_{i,j}v_{i,k}) + \sum^{n+1}_{k=j+1} f(v_{i,j}v_{i,k})\nonumber\\
&=& m{n+1\choose2}+(m-i+1)n - j+1   +  \sum^{j-1}_{k=1} \Bigg\{m\Bigg[{j-1\choose2}+k-1\Bigg]+i\Bigg\}  + \nonumber\\
& & \sum^{n+1}_{k=j+1} \Bigg\{m\Bigg[{k-1\choose2}+j-1\Bigg]+i\Bigg\}  \nonumber\\
&=& m{n+1\choose2}+(m+1)n - j+1 +  mj{j-1\choose2} + m{n+1\choose3}-m{j\choose3}+\nonumber\\
& & m(n-j+1)(j-1)\nonumber\\
&=& m{n+1\choose2}+(m+1)n -j + 1 + m{n+1\choose3} + mn(j-1) + \nonumber\\ 
& & m\Bigg[j^3/3-2j^2+8j/3-1\Bigg]
\end{eqnarray} 
also an increasing function.

\nt Thus, we have $w(v_{i,1}) < w(v_{i,2}) <\cdots < w(v_{i,n})$ and $w(v_{1,n+1}) < w(v_{2,n+1}) < \cdots < w(v_{m,n+1})$. It is routine to check that $$w(v_{1,n+1}) - w(v_{i,n}) = \frac{mn}{2}(mn+n+3m-7)+n+\frac{m}{6}(m^2+3m+8)-1 > 0.$$ Therefore, $f$ is a local antimagic total labeling that induces $m+n$ distinct weights so that $\chi_{lat}(K_m\odot K_n)\le m+n$. Since $\chi(K_m \odot K_n)=\max\{m,n+1\}$, we have $\max\{m,n+1\}\le \chi_{lat}(K_m\odot K_n)\le m+n$. When $m=1$, we have $\max\{m,n+1\}=n+1$. The theorem holds.
\end{proof}

\nt In what follows, we consider some sparse graphs. Let $P_n = v_1v_2\cdots v_n$ be the path of order $n\ge 2$. 

\begin{theorem}\label{thm-Pn} For $n\ge 2$, $\chi_{lat}(P_n)=2$ except that $\chi_{lat}(P_4)=3$. 
\end{theorem}

\begin{proof} Obviously $\chi_{lat}(P_2)=2$. By~\cite[Theorem 3.5]{LSN1}, $\chi_{la}(P_4\vee K_1) = 4$ implies that $\chi_{lat}(P_4)\ne 2$. Moreover, every local antimagic labeling of $P_4\vee K_1$ induces 4 distinct vertex colors also corresponds to a local antimagic total labeling of $P_4$ that induces 3 distinct vertex weights. One such labeling that labeled the vertices and edges  of $P_4$ alternately is given by sequence $7$, $2$, $6$, $3$, $5$, $1$, $4$.

\ms\nt Assume $n\ge 6$ is even. Observe that by Theorem~\ref{thm-K1VG} and the local antimagic labeling of $W_n$ obtained in~\cite{Arumugam,LNS}, we can get a local antimagic total labeling of $C_n$ with an edge labeled $1$. Since $C_n$ is regular, we can delete this edge and reduce all other labels by 1. Consequently, we get a path $P_n$ that admits a local antimagic total labeling that induces exactly 2 distinct vertex weights. Since $\chi_{lat}(P_n)\ge \chi(P_n)\ge 2$, we have $\chi_{lat}(P_n)=2$. The details are omitted. We note that similar conclusion can also be obtained from the labeling of $W_n, n$ even in Theorem~\ref{thm-A(mWn,v)} for $m=1$ after relabeling the vertices and edges of $W_n$ in reverse order. 

\ms\nt Consider odd $n\ge 5$. Suppose $n=4k+1$. For $n=5$, a required labeling sequence that labeled the vertices and edges of $P_5$ alternately is $6$, $4$, $7$, $3$, $2$, $5$, $8$, $1$, $9$  with distinct vertex weights 10 and 14. By computer search, we are able to obtain all the 12 different labelings. For $n\ge 9$, define $f : V(P_n) \cup E(P_n) \to [1, 8k+1]$ as follows.

\begin{enumerate}[(i)]
  \item $f(v_{1}) = 8k$; $f(v_{4k-1}) = 6k$; $f(v_{4k+1}) = 8k+1$; 
  \item $f(v_{4i+1}) = 7k-3+i$ for $i\in [1,k-1]$;
  \item $f(v_{4i-1}) = 3k+i$ for $i\in [1,k-1]$;
  \item $f(v_{4i+2}) = 7k+i$ for $i\in [0,k-1]$;
  \item $f(v_{4i+4}) = 4k+1+i$ for $i\in [0,k-1]$;
  \item $f(v_{4k}v_{4k+1}) = 2k$; $f(v_{4k-1}v_{4k}) = 4k$; 
  \item $f(v_{2i}v_{2i+1}) = 2k-i$ for $i\in [1,2k-1]$;
  \item $f(v_{4i+1}v_{4i+2}) = 2k+1+i$ for $i\in [0,k-1]$;
  \item $f(v_{4i+3}v_{4i+4}) = 5k+1+i$ for $i\in [0,k-2]$.
\end{enumerate}

\ms\nt It is not difficult to check that $$w(v_i) =\begin{cases} 10k+1 & \mbox{ for odd } i,\\ 11k & \mbox{ for even } i.\end{cases}$$ Thus, $\chi_{lat}(P_{4k+1}) = 2$. 

\ms\nt Suppose $n=4k+3\ge 3$. A required labeling sequence for $n=3,7,11$ are $1$, $5$, $3$, $4$, $2$; $5$, $13$, $1$, $9$, $7$, $2$, $10$, $11$, $4$, $3$, $8$, $12$, $6$; and $8$, $21$, $1$, $14$, $11$, $4$, $15$, $17$, $10$, $2$, $16$, $18$, $6$, $5$, $12$, $19$, $7$, $3$, $12$, $20$, $9$. The corresponding distinct vertex weights are 6 and 12; 18 and 23; and 29 and 36 respectively. 

\ms\nt For $k\ge 3$, we define $f: V(P_n)\cup E(P_n)\to [1, 8k+5]$ as follows.

\begin{enumerate}[(i)]
  \item $f(v_1) = 3k+2$; $f(v_2) = 1$; $f(v_3) = 4k+3$; $f(v_{2k+3}) = 2k+2$; $f(v_{4k+3}) = 3k+3$;
  \item $f(v_{2i+3}) = 3k+3+i$ for $i\in[1,k-1]$;
  \item $f(v_{2k+2i+1}) = 2k+1+i$ for $i\in[1,k]$;
  \item $f(v_{2i+2}) = 5k+4+i$ for $i\in[1,k]$;
  \item $f(v_{2k+2i+2}) = 4k+3+i$ for $i\in[1,k]$;
  \item $f(v_1v_2) = 8k+5$; $f(v_2v_3) = 5k+4$;
  \item $f(v_{2i+1}v_{2i+2}) = 2k+2-2i$ for $i\in[1,k]$;
  \item $f(v_{2k+2i+1}v_{2k+2i+2}) = 2k+3-2i$ for $i\in[1,k]$;
  \item $f(v_{2i+2}v_{2i+3}) = 6k+4+i$ for $i\in[1,2k]$.
\end{enumerate}

\ms\nt It is not difficult to check that $$w(v_i) =\begin{cases} 11k+7 & \mbox{ for odd } i,\\ 13k+10 & \mbox{ for even } i.\end{cases}$$ Thus, $\chi_{lat}(P_{4k+3}) = 2$.
\end{proof}

\begin{example}\label{ex-path} The labeling sequence for $P_{12}$ is $16$, $10$, $18$, $3$, $12$, $11$, $19$, $1$, $20$, $5$, $17$, $9$, $15$, $2$, $22$, $7$, $13$, $6$, $21$, $4$, $14$, $8$, $23$ with $2$ distinct vertex weights $26$ and $31$. The labeling sequence for $P_{16}$ is $22$, $13$, $24$, $5$, $16$, $14$, $25$, $3$, $17$, $15$, $26$, $1$, $27$, $7$, $23$, $12$, $21$, $2$, $30$, $10$, $19$, $6$, $28$, $8$, $18$, $9$, $29$, $4$, $20$, $11$, $31$ with $2$ distinct vertex weights $35$ and $42$.

\ms\nt The labeling sequence for $P_{14}$ is $24$, $13$, $16$, $1$, $27$, $9$, $19$, $2$, $23$, $12$, $15$, $3$, $26$, $8$, $18$, $4$, $22$, $11$, $14$, $5$, $25$, $7$, $17$, $6$, $21$, $10$, $20$ with $2$ distinct vertex weights $37$ and $30$.

\ms\nt The labeling sequence for $P_{13}$ is $24$, $7$, $21$, $5$, $10$, $16$, $13$, $4$, $19$, $8$, $22$, $3$, $11$, $17$, $14$, $2$, $20$, $9$, $23$, $1$, $18$, $12$, $15$, $6$, $25$ with $2$ distinct vertex weights $31$ and $33$.

\ms\nt The labeling sequence for $P_{15}$ is $11$, $29$, $1$, $19$, $15$, $6$, $20$, $23$, $13$, $4$, $21$, $24$, $14$, $2$, $22$, $25$, $8$, $7$, $16$, $26$, $9$, $5$, $17$, $27$, $10$, $3$, $18$, $28$, $12$ with $2$ distinct vertex weights $40$ and $49$. \end{example}

\ms\nt In~\cite[Theorems 3.6 and 3.7]{LSN1}, the authors showed that $\chi_{la}(F_n)=3$ for even $n\ge 4$, and $3\le\chi_{la}(F_n)\le 4$ for odd $n\ge 3$. From the proof of Theorem~\ref{thm-Pn}, when $n=4k+1$, we have $\sum^{4k+1}_{i=1} f(u_i) = 22k^2 + 12k + 1 \ne 11k \ne 10k+1$. For $n=4k+3$, we have $\sum^{4k+3}_{i=1} f(u_i) = 16k^2 + 19k + 6 \ne 11k+7 \ne 13k+10$. By Theorems~\ref{thm-Pn} and~\ref{thm-K1VG}, we have the following.

\begin{corollary} For $n\ge 2$, $\chi_{la}(F_{n}) = 3$. 
\end{corollary}


\begin{theorem}\label{thm-Cn} For $n\ge 3$, $$\chi_{lat}(C_n)=\begin{cases}2 &\mbox{ if } n \mbox{ is even,}\\ 3 &\mbox{ otherwise.}\end{cases}$$
\end{theorem}

\begin{proof} It is obvious that $\chi_{lat}(C_3) = 3$. Assume $n\ge 4$. In~\cite{Arumugam, LNS}, the authors showed that $$\chi_{la}(W_n)=\begin{cases}3 &\mbox{ if } n \mbox{ is even,}\\ 4 &\mbox{ otherwise.}\end{cases}$$ Since $$\chi(C_n)=\begin{cases}2 &\mbox{ if } n \mbox{ is even,}\\ 3 &\mbox{ otherwise,}\end{cases}$$ by Theorem~\ref{thm-K1VG}, we conclude that the theorem holds.
\end{proof}





\begin{theorem} For odd $n\ge 3$, $4\le \chi_{lat}(C_n \vee 2K_1)\le 5$, and for even $n\ge 6$, $3\le \chi_{lat}(C_n \vee 3K_1)\le 5$. \end{theorem}

\begin{proof} Suppose $n\ge 3$ is odd. Clearly, $\chi_{lat}(C_n \vee 2K_1)\ge \chi(C_n \vee 2K_1)=4$. In~\cite[Theorem 3.1]{LSN1}, the authors proved that for odd $n\ge 3$, $\chi_{la}(C_n\vee 3K_1) = 4$. Moreover, the corresponding local antimagic labeling $g$ induces $g^+(u_1)=g^+(u_2)=g^+(u_3)=5n(n+1)/2$, $g^+(v_1)=8n+3$, $g^+(v_i)= (17n+7)/2$ for odd $i\ge 3$, and $g^+(v_i) = (17n+5)/2$ for even $i\ge 2$.

\ms\nt Define $f:V(C_n\vee 2K_1) \cup E(C_n\vee 2K_1) \to [1,4n+2]$ such that $f(v_i)=g(v_iu_3)$, and $f(e) = g(e)$ for $e\in E(C_n)$ or $e=v_iu_j, j=1,2$. Moreover, $f(u_j)=4n+j$ for $j=1,2$. Now, $w(v_i)=g^+(v_i)$ and $w(u_j)=g^+(u_j)+4n+i$ for $i=1,2$. Thus, $f$ induces 5 distinct vertex weights and $\chi_{lat}(C_n \vee 2K_1)\le 5$. 

\ms\nt Suppose $n\ge 6$ is even. Clearly, $\chi_{lat}(C_n \vee 3K_1)\ge 3$. In~\cite[Theorem 3.3]{LSN1}, the authors proved that $\chi_{la}(C_n\vee 4K_1)=3$. Moreover, the corresponding local antimagic labeling $g$ induces $g^+(v_i) = 9n+3$ for odd $i$, $g^+(v_i) = 17n+3$ for even $i$, and $g^+(u_j)=n(6n+1)/2$ for $1\le j\le 4$.  

\ms\nt Define $f:V(C_n \vee 3K_1)\cup E(C_n \vee 3K_1)\to [1,5n+3]$ such that $f(v_i) = g(v_iu_4)$, and $f(e) = g(e)$ for $e\in E(C_n)$ or $e=v_iu_j, j=1,2,3$. Moreover, $f(u_j)=5n+j$ for $j=1,2,3$. Now $w(v_i) = g^+(v_i)$ and $w(u_j)=g^+(u_j)+5n+i$ for $i=1,2,3$. Thus, $f$ induces 5 distinct vertex weights and  $\chi_{lat}(C_n \vee 3K_1)\le 5$.  
\end{proof} 

\begin{problem} Determine $\chi_{lat}(C_n \vee 2K_1)$ for odd $n\ge 3$, and $\chi_{lat}(C_n \vee 3K_1)$ for even $n\ge 4$. \end{problem}

\ms\nt In~\cite[Theorem 3.9]{LSN1}, the authors proved that for $n,m\ge 3$, $$\chi_{la}(K_m\vee C_n) = \begin{cases}m+2 & \mbox{ if } m,n  \mbox{ are even;} \\ m+3 &\mbox{ if } m,n \mbox{ are odd.}\end{cases}$$ By Theorem~\ref{thm-K1VG}, the following theorem holds. 

\begin{theorem}
For $m,n\ge 3$, $$\chi_{lat}(K_{m-1}\vee C_n) = \begin{cases} m+1 & \mbox{ if } m,n \mbox{ are even;} \\ m+2 &\mbox{ if } m,n \mbox{ are odd.}\end{cases}$$
\end{theorem}

\nt In~\cite[Theorem 4]{Slamin+NMDK}, the authors also proved that a family of cubic bipartite graph denoted $CB_{2k}$ obtained from $C_{2k}, k\ge 3$ has $\chi_{lsat}(CB_{2k})=2$. Thus, by Lemma~\ref{lem-regular}, we have 

\begin{corollary} For $k\ge 3$, $\chi_{lat}(CB_{2k})=\chi_{lat}(CB_{2k}-e) = 2$. \end{corollary}

\nt It is easy to verify that the given labeling has sum of all the vertex labels not equal the induced vertex weights. By Theorem~\ref{thm-K1VG} (b), we have

\begin{corollary} For $k\ge 3$, $\chi_{la}(CB_{2k} \vee K_1) = 3$. \end{corollary}

\nt By Lemma~2.4 in~\cite{LSN1}, it is easy to check that the following corollary holds.

\begin{corollary} For $k\ge 3$, $\chi_{la}((CB_{2k} \vee K_1)-e)=3$ where $e$ is not an edge that belongs to the induced $C_{2k}$ subgraph. \end{corollary}

\nt Note that when $k$ is odd, $CB_{2k}$ is commonly known as M\"obius ladder $M_{2k}$ or the circulant graph $C_{2k}(1,k)$.

\section{Conclusion and Open Problems}

\nt In this paper, we first prove that every graph is local antimagic. The proof gives a good bound for us to determine $\chi_{lat}(G)$ (or $\chi_{la}(G\vee K_1)$) using a local antimagic labeling of $G\vee K_1$ (or a local antimagic total labeling of $G$). The local antimagic (total) chromatic number of many family of graphs are determined. Particularly, we showed that there are graphs $G$ with $\chi(G)=\chi(G-e)=\chi_{lat}(G)=\chi_{lat}(G-e)=\chi(G\vee K_1)-1=\chi_{la}(G\vee K_1)-1=\chi_{la}((G-e)\vee K_1)-1$ or $\chi(G)=\chi_{lat}(G) = \chi_{la}(G)-1$. The following problems arise naturally.

\begin{problem} Characterize $G$ such that $\chi(G)=\chi(G-e)=\chi_{lat}(G)=\chi_{lat}(G-e)=\chi(G\vee K_1)-1=\chi_{la}(G\vee K_1)-1=\chi_{la}((G-e)\vee K_1)-1$. \end{problem}

\begin{problem} Characterize $G$ such that $\chi(G)=\chi_{lat}(G) = \chi_{la}(G)-1$. \end{problem}

\begin{problem} Determine the exact values of $\chi_{lat}(A(mK_n,K_r)$ and $\chi_{lat}(A(mK_n,K_r)-e)$ for $m,n,r$ not satisfying Theorems~\ref{thm-A(mKnKr)} and~\ref{thm-A(mKnKr-e)}.\end{problem}

\begin{problem} For $n\ge 3$, prove that if $e$ is an edge adjacent to the vertex $u$ of $A(2K_n,K_1)$, then $\chi_{lat}(A(2K_n,K_1)-e)=n$. \end{problem}

\begin{problem} Determine the exact values of $\chi_{lat}(mK_n)$ and $\chi_{lat}(mK_n - e)$ for $m\ge 2$ and odd $n\ge 3$. \end{problem}

\begin{problem} Determine the exact values of $\chi_{lat}(K_m\odot K_n)$ and $\chi_{la}((K_m\odot K_n)\vee K_1)$ for $m,n\ge 2$. \end{problem}

\begin{problem} Determine the exact values of $\chi_{lat}(F_n)$ and $\chi_{la}(F_n \vee K_1)$ for odd $n\ge 3$. \end{problem}

\begin{problem} Determine the exact values of $\chi_{lat}(C_n\vee 2K_1)$ and $\chi_{la}((C_n\vee 2K_1)\vee K_1)$ for odd $n\ge 3$, and $\chi_{lat}(C_n\vee 3K_1)$ and $\chi_{la}((C_n\vee 3K_1)\vee K_1)$ for even $n\ge 6$. \end{problem}

\nt In~\cite[Theorem 3.4]{Lau+LNS}, the authors showed that there are infinitely many circulant graphs (with at most an edge deleted) of $\chi_{la}=3$. Since cycles are the simplest circulant graphs with $\chi_{lat}=2$, we have

\begin{problem} Determine the exact values of $\chi_{lat}(C)$ and $\chi_{lat}(C-e)$ for each circulant graph $C\not\cong C_n, C_{2n}(1,n), n\ge 3$. \end{problem}

\nt In~\cite{Lau+SN}, the authors proved that any graph $G$ of $k\ge 1$ pendant(s) has $\chi_{la}(G)\ge k+1$. In~\cite{Arumugam, Lau+LNS, LNS, LSN1, Lau+SS, Premalatha+ALW}, the authors obtained many families of graph $G$ having $k\ge 1$ pendants with $\chi_{la}(G)=k+i, i=1,2$. It is also obvious that $\chi_{la}(C_3 \odot O_2)=9$. Besides $\chi_{la}(K_n)=n, n\ge 3$, we are not aware of any graph $G\not\cong K_n$ containing no pendant vertices but having arbitrarily large $\chi_{la}(G)=\chi(G)$. Our results above show the existence of infinitely many such graphs $G$. This gives partial solution to \cite[Problem 3.2]{Arumugam}: Characterize the class of graphs $G$ for which $\chi_{la}(G)=\chi(G)$. 

\ms\nt Since every known result has $\chi_{lat}(G)\le\chi_{la}(G)$, we end this paper with the following.

\begin{conjecture} For each graph $G$ of order at least 3, $\chi_{lat}(G)\le \chi_{la}(G)$. \end{conjecture}


\begin{thebibliography}{99}
\bibitem{Arumugam} S. Arumugam, K. Premalatha, M. Bac\v{a} and A. Semani\v{c}ov\'{a}-Fe\v{n}ov\v{c}\'{i}kov\'{a}, Local antimagic vertex coloring of a graph, {\it Graphs and Combin.}, {\bf33} (2017), 275-285.

\bibitem{Bensmail} J. Bensmail, M. Senhaji and K. Szabo Lyngsie, On a combination of the 1-2-3 Conjecture and the Antimagic Labelling Conjecture, {\it Discrete Math. Theoret. Comput. Sc.}, {\bf19(1)} (2017) \#22.



\bibitem{Haslegrave} J. Haslegrave, Proof of a local antimagic conjecture, {\it Discrete Math. Theor. Comp. Sc.}, {\bf 20(1)} (2018), \#18.

\bibitem{Lau+LNS} G.C. Lau, J. Li, H.K. Ng and W.C. Shiu, Approaches which output infinitely many graphs with small local antimagic chromatic number, (2020), submitted to {\it Discrete Math}, arXiv:2009.01996.

\bibitem{LNS} G.C. Lau, H.K. Ng, and W.C. Shiu, Affirmative solutions on local antimagic chromatic number, {\it Graphs and Combin.}, {\bf36} (2020), 1337-1354.


\bibitem{Lau+SN} G.C. Lau, W.C. Shiu and H.K. Ng, On local antimagic chromatic number of graphs with cut-vertices, (2020) submitted to {\it Iran. J. Math. Sci. Inform.}, arXiv:1805.04801.

\bibitem{LSN1} G.C. Lau, W.C. Shiu and H.K. Ng, On local antimagic chromatic number of cycle-related join graphs, {\it Discuss. Math. Graph Theory} {\bf 41} (2021) 133–152 DOI : 10.7151/dmgt.2177.

\bibitem{Lau+SS} G.C. Lau, W.C. Shiu and C.X. Soo, On Local Antimagic Chromatic Number of Spider Graphs, (2020) submitted to {\it J. Discret. Math. Sci. Cryptogr.},  arXiv:2008.09754.

\bibitem{Premalatha+ALW} K. Premalatha, S. Arumugam, Yi-Chun Lee and Tao-Ming Wang, Local antimagic chromatic number of trees - I, {\it J. Discret. Math. Sci. Cryptogr.}, (2020) DOI : 10.1080/09720529.2020.1772985. 

\bibitem{Slamin+NMDK} Slamin, N. O. Adiwijaya, M. A. Hasan, Dafik, and K. Wijaya, Local Super Antimagic Total Labeling for Vertex Coloring of Graphs, {\it Symmetry}, {\bf 12(11)} (2020) 10.3390/sym12111843.

\bibitem{Zukerman} D. Zuckerman, Linear degree extractors and the inapproximability of Max Clique and Chromatic Number, {\it Theory of Computing}, {\bf3} (2007) 103-128, doi:10.4086/toc.2007.v003a006
\end{thebibliography}
\end{document}